\nonstopmode
\documentclass[10pt]{amsart}
\usepackage{latexsym}
\usepackage{fancyhdr}
\usepackage{amsmath, amssymb}
\usepackage[ansinew]{inputenc}
\usepackage[all]{xy}
\usepackage{pdflscape}
\usepackage{longtable}
\usepackage{rotating}
\usepackage{verbatim}
\usepackage{hyperref}
\usepackage{subfigure}
\usepackage{mathrsfs}
\usepackage{enumerate}
\usepackage{mdwlist}

\swapnumbers
\theoremstyle{plain}
\newtheorem*{lemma*}{Lemma}
\newtheorem{lemma}[subsection]{Lemma}
\newtheorem*{theorem*}{Theorem}
\newtheorem{theorem}[subsection]{Theorem}
\newtheorem*{proposition*}{Proposition}
\newtheorem{proposition}[subsection]{Proposition}
\newtheorem*{corollary*}{Corollary}
\newtheorem{corollary}[subsection]{Corollary}
\newtheorem*{claim*}{Claim}
\newtheorem{claim}[subsection]{Claim}

\newtheorem*{conjecture*}{Conjecture}

\newtheorem*{question*}{Question}
\theoremstyle{definition}
\newtheorem*{definition*}{Definition}

\newtheorem*{example*}{Example}
\newtheorem{example}[subsection]{Example}
\newtheorem{examples}[subsection]{Examples}

\newtheorem*{algorithm*}{Algorithm}
\newtheorem*{remark*}{Remark}
\newtheorem*{remarks*}{Remarks}
\newtheorem{remark}[subsection]{Remark}
\newtheorem{remarks}[subsection]{Remarks}
\newtheorem*{convention*}{Convention}

\makeatletter
  \let\c@equation\c@subsection
  
\makeatother

\newenvironment{demo}[1]{\par\smallskip\noindent{\bf #1.}}{\par\smallskip}
\numberwithin{equation}{section}

\def\al{\alpha}
\def\be{\beta}
\def\ga{\gamma}

\def\la{\lambda}
\def\vp{\varpi}
\def\rh{\rho}

\def\si{\sigma}

\def\ta{\tau}

\def\vh{\varphi}

\def\om{\omega}

\def\Ga{\Gamma}

\def\La{\Lambda}
\def\Si{\Sigma}

\def\Om{\Omega}

\def\C{\mathbb{C}}

\def\N{\mathbb{N}}

\def\R{\mathbb{R}}

\def\cB{\mathcal{B}}

\def\cE{\mathcal{E}}
\def\cF{\mathcal{F}}
\def\cG{\mathcal{G}}
\def\cH{\mathcal{H}}

\def\fG{\mathfrak{G}}

\def\fM{\mathfrak{M}}
\def\fN{\mathfrak{N}}

\def\fS{\mathfrak{S}}

\def\fW{\mathfrak{W}}

\def\sM{\mathscr{M}}

\def\sW{\mathscr{W}}

\def\p{\partial}

\def\<{\langle}
\def\>{\rangle}
\renewcommand{\o}{\circ}

\let\on=\operatorname

\newcommand{\sr}[1]%
{\ifmmode{}^\dagger\else${}^\dagger$\fi\ifvmode
\vbox to 0pt{\vss
 \hbox to 0pt{\hskip\hsize\hskip1em
 \vbox{\hsize3cm\raggedright\pretolerance10000
 \noindent #1\hfill}\hss}\vss}\else
 \vadjust{\vbox to0pt{\vss%
 \hbox to 0pt{\hskip\hsize\hskip1em%
 \vbox{\hsize3cm\raggedright\pretolerance10000%
 \noindent #1\hfill}\hss}\vss}}\fi%
}

\providecommand{\mapsfrom}{\kern.2em%
\setbox0=\hbox{$\leftarrow$\kern-.10em\rule[0.26mm]{0.1mm}{1.3mm}}\box0%
\kern.3em}

\title[Composition in ultradifferentiable classes]
{Composition in ultradifferentiable classes}

\author[A.~Rainer]{Armin Rainer}

\address{A.~Rainer: 
Fakult\"at f\"ur Mathematik, Universit\"at Wien, 
Oskar-Morgenstern-Platz~1, A-1090 Wien, Austria}
\email{armin.rainer@univie.ac.at}

\author[G.~Schindl]{Gerhard Schindl}

\address{G.~Schindl: Fakult\"at f\"ur Mathematik, Universit\"at Wien, 
Oskar-Morgenstern-Platz~1, A-1090 Wien, Austria}
\email{a0304518@unet.univie.ac.at}

\begin{document}

\begin{abstract}
  We characterize stability under composition of ultradifferentiable classes defined by 
  weight sequences $M$, by weight functions $\om$, and, more generally, by weight matrices $\fM$,
  and investigate continuity of composition $(g,f) \mapsto f \o g$. 
  In addition, we represent the Beurling space $\cE^{(\om)}$ and the Roumieu space $\cE^{\{\om\}}$ as 
  intersection and union of spaces $\cE^{(M)}$ and $\cE^{\{M\}}$ for associated weight sequences,
  respectively.   
\end{abstract}


\thanks{AR was supported by FWF-Project P~22218-N13 and GS by FWF-Project P~23028-N13}
\date{November 24, 2013}

\maketitle

\section{Introduction}

This paper arose from our wish to characterize stability under composition of Denjoy--Carleman classes 
$\cE^{\{M\}}$ and $\cE^{(M)}$. For these classes we have developed a calculus in infinite dimensions beyond Banach 
spaces in \cite{KMRc,KMRq,KMRu} which is heavily based on composition: 
A smooth mapping $f$ is of class $\mathcal E^{\{M\}}$ if and only if $f \circ p \in \mathcal E^{\{M\}}$ 
for all  $\mathcal E^{\{M\}}$ Banach plots (i.e., mappings defined in open subsets of Banach spaces); 
accordingly for  $\mathcal E^{(M)}$. Sometimes curves suffice.

Denjoy--Carleman differentiable functions form classes of smooth functions that are 
described by growth conditions on the Taylor expansion.
The growth is prescribed in terms of a sequence $M=(M_k)$ of positive real numbers which 
serves as a weight for the iterated derivatives: for compact $K$ the sets 
\[
\Big\{\frac{f^{(k)}(x)}{\rh^{k} \, k! \, M_{k}} : x \in K, k \in \N \Big\}
\] 
are required to be bounded. The positive real number $\rh$ is subject to either 
a universal or an existential quantifier, thereby dividing the Denjoy--Carleman classes into 
those of Beurling type $\cE^{(M)}$ and those of Roumieu type $\cE^{\{M\}}$, 
respectively. 
We write $\cE^{[M]}$ for either $\cE^{(M)}$ or $\cE^{\{M\}}$.

It is well-known that $\cE^{[M]}$ is stable under composition, if $M$ is log-convex, 
see \cite{Roumieu62/63}, \cite{Komatsu73}, \cite{Dynkin80}, and usually in the literature log-convexity is assumed in 
order to have stability under composition; but is log-convexity also necessary? 
Actually, when proving stability under composition with Fa\'a di Bruno's 
formula one needs a weaker condition that we call (FdB)-property. 
We prove that the (FdB)-property (for the weakly log-convex minorant $M^{\flat(c)}$) is also a necessary 
condition for stability under composition, if $\cE^{[M]}$ 
is stable under derivation, 
see Theorem~\ref{thm:comp_dc}. More precisely, if $\cE^{[M]}$ 
is stable under derivation, then  
stability under composition is in turn equivalent to 
being holomorphically closed, 
being inverse closed, $(M^{\flat(c)}_k)^{\frac{1}{k}}$ being almost increasing, and $M^{\flat(c)}$ having the (FdB)-property. 
For further equivalent stability properties we refer to \cite{RainerSchindl14}.
Inverse closedness has been studied intensively, e.g.\ \cite{Rudin62}, \cite{Bruna80/81}, \cite{Siddiqi90}.
In this context we prove that, as in the Roumieu case \cite{Cartan40}, one has $\cE^{(M)}=\cE^{(M^{\flat(c)})}$ if 
$C^\om \subseteq \cE^{(M)}$, see Theorem~\ref{thm:reg}.  
Finally, we demonstrate that log-convexity is not necessary for stability under composition: 
We construct classes $\cE^{[M]}$ which are stable under composition and such that there is no log-convex 
$N=(N_k)$ with $\cE^{[M]} = \cE^{[N]}$, see Example~\ref{ex:lc}.   
       
Another common way to define ultradifferentiable classes is by means of a weight function $\om$ which controls the decay of the Fourier
transform, see \cite{Beurling61} and \cite{Bjoerck66}; we shall use the following equivalent description due to \cite{BMT90}: 
for compact $K$ the sets 
\[
\Big\{f^{(k)}(x) \exp(-\tfrac{1}{\rh} \vh^*(\rh k)) : x \in K, k \in \N \Big\},
\] 
where $\vh^*$ is the Young conjugate of $\vh(t) = \om(e^t)$,
are required to be bounded either for all $\rh>0$ in the Beurling case $\cE^{(\om)}$ 
or for some $\rh>0$ in the Roumieu case $\cE^{\{\om\}}$. 
Again $\cE^{[\om]}$ stands for either $\cE^{(\om)}$ or $\cE^{\{\om\}}$.
For these classes stability under composition was characterized in \cite{FernandezGalbis06} under the additional assumption of 
non-quasianalyticity. 
Note that the sets $\{\cE^{[M]} : M \text{ weight sequence}\}$ and $\{\cE^{[\om]} : \om \text{ weight function}\}$ 
have a \emph{large} intersection but neither of them contains the other, see 
\cite{BMM07}. 
We want to stress the fact that the usual requirements on the weight function $\om$ ensure that the spaces $\cE^{[\om]}$ 
come with incorporated stability properties, for instance stability under derivation, see Corollary \ref{cor:rep}. 

We prove that $\cE^{(\om)}$ and $\cE^{\{\om\}}$ can be represented (as locally convex spaces with their natural topologies) 
as intersections and unions of ultradifferentiable classes defined by means of associated weight sequences, see Theorem~\ref{thm:rep}:
For each open subset $U \subseteq \R^n$, compact $K \subseteq U$, and for $\Om^\rh = (\Om^\rh_k)$ defined by 
$\Om^\rh_k:= \tfrac{1}{k!} \exp(\tfrac{1}{\rh} \vh^*(\rh k))$ we have   
\begin{equation} \label{eq:rep}
  \cE^{(\om)}(U) = \bigcap_{\rh>0} \cE^{(\Om^\rh)}(U) \quad \text{ and } \quad  
  \cE^{\{\om\}}(U) = \bigcap_{K \subseteq U} \bigcup_{\rh>0} \cE^{\{\Om^\rh\}}(K).  
\end{equation}
We use this representation for characterizing stability under composition, and believe that it is also of independent interest.

In fact, inspired by \eqref{eq:rep}, we characterize stability under composition for more general ultradifferentiable classes 
defined by weight matrices $\fM = \{M^\la \in \R_{>0}^\N : \la \in \La\}$, where 
$\La$ is 
an ordered subset of $\R$:
\begin{equation*} 
  \cE^{(\fM)}(U) := \bigcap_{\la \in \La} \cE^{(M^\la)}(U) \quad \text{ and } \quad  
  \cE^{\{\fM\}}(U) := \bigcap_{K \subseteq U} \bigcup_{\la \in \La} \cE^{\{M^\la\}}(K),  
\end{equation*}
endowed with their natural topologies. Among the spaces $\cE^{(\fM)}$ and $\cE^{\{\fM\}}$,
commonly denoted by $\cE^{[\fM]}$, are all 
the spaces defined by means of weight sequences and weight functions, but not exclusively, see Theorem~\ref{thm:ext}. 
For instance, the intersection, resp.\ the union, of all non-quasianalytic Gevrey classes is an autonomous 
$\cE^{(\fM)}$-space, resp.\ $\cE^{\{\fM\}}$-space, with suitable $\fM$.
Intersections of non-quasianalytic ultradifferentiable classes have been studied by 
Rudin \cite{Rudin62}, Boman \cite{Boman63}, 
Chaumat and Chollet \cite{ChaumatChollet98}, 
Beaugendre \cite{Beaugendre01, Beaugendre02}, and Schmets and Valdivia \cite{SchmetsValdivia05/06,SchmetsValdivia08}
(among others). 
It seems, however, that unions of ultradifferentiable classes have not been investigated before.

Given that $\cE^{[\fM]}$ is stable under composition, the nonlinear composition operators 
\begin{gather*}
 \on{comp}^{(\fM)} : \cE^{(\fM)}(\R^p,\R^q) \times \cE^{(\fM)}(\R^q,\R^r) \to \cE^{(\fM)}(\R^p,\R^r) : (g,f) \mapsto f \o g \\
 \cE^{\{\fM\}}(\R^p,f) : \cE^{\{\fM\}}(\R^p,\R^q) \to \cE^{\{\fM\}}(\R^p,\R^r) : g \mapsto f \o g, \quad f \in \cE^{\{\fM\}}(\R^q,\R^r),
\end{gather*}
turn out to be continuous. This is proved in Theorem~\ref{thm:compC0}. The special case $\cE^{[\om]}$ was treated in 
\cite{FernandezGalbis06}, see also \cite{AppellNazarovZabrejko91}. 
Under suitable assumptions we expect $\on{comp}^{[\fM]}$ to be of class $\cE^{[\fM]}$, see Remark~\ref{rem:fc}.

The paper is structured as follows:
We first treat the weight sequence case in Section~\ref{sec:ws} and Section~\ref{sec:wsc}.
In Section~\ref{sec:wm} we introduce ultradifferentiable classes defined by weight matrices $\fM$, 
characterize their stability under composition, and show that composition is continuous. 
We discuss classes defined by weight functions $\om$ and identify them as classes defined by weight matrices $\fM$ in Section~\ref{sec:wf},
and characterize their stability under composition in 
Section~\ref{sec:wfc}.

\subsection*{Notation and conventions}
The notation $\cE^{[*]}$ for $* \in \{M,\om,\fM\}$ stands for either $\cE^{(*)}$ or $\cE^{\{*\}}$ with the following restriction: 
Statements that involve more than one $\cE^{[*]}$ symbol must not be interpreted by mixing $\cE^{(*)}$ and $\cE^{\{*\}}$. 
This convention will be used broadly, but self-evidently: 
For example, $\fM [\preceq] \fN \Leftrightarrow \cE^{[\fM]} \subseteq \cE^{[\fN]}$ 
in Proposition~\ref{prop:fMincl} means $\fM (\preceq) \fN \Leftrightarrow \cE^{(\fM)} \subseteq \cE^{(\fN)}$ 
and $\fM \{\preceq\} \fN \Leftrightarrow \cE^{\{\fM\}} \subseteq \cE^{\{\fN\}}$.   

Let $\N = \N_{>0} \cup \{0\}$. 
For $\al=(\al_1,\ldots,\al_q) \in \N^q$ and $x = (x_1,\ldots,x_q) \in \R^q$
we write $\al!=\al_1! \cdots \al_q!$, $|\al|= \al_1 +\cdots+ \al_q$, and $x^\al = x_1^{\al_1}\cdots x_q^{\al_q}$.
We use $\p_i = \p/\p x_i$, 
$\p^\al= \p_1^{\al_1}\cdots \p_q^{\al_q}$ and write 
$d^k f$ or 
$f^{(k)}$ for the $k$th order Fr\'echet derivative of $f$, and $d_v f$ for the directional derivative in direction $v$.
For sequences of reals $M=(M_k)$ and $N=(N_k)$ we write $M \le N$ if $M_k \le N_k$ for all $k$.

$L(E_1,\ldots,E_k;F)$ is the space of $k$-linear bounded mappings 
$E_1 \times \cdots \times E_k \to F$ (between topological vector spaces); if $E_i =E$ for all $i$, we also write $L^k(E,F)$.

Let $\cF$ and $\cG$ denote classes of mappings.  
We write $\cF \subseteq \cG$ if $\cF(U,\R^m) \subseteq \cG(U,\R^m)$ for all open subsets $U \subseteq \R^n$ and all $n,m \in \N_{>0}$. 
We say that \emph{$\cF$ is stable under composition} 
if $g \in \cF(U,V)$ and $f \in \cF(V,W)$ implies $f\o g \in \cF(U,W)$, 
for all open subsets $U \subseteq \R^p$, $V \subseteq \R^q$, $W \subseteq \R^r$, and all $p,q,r \in \N_{>0}$.
A class $\cF$ is called \emph{holomorphically closed} if $f \o g \in \cF(U,\C)$ for each $g \in \cF(U) = \cF(U,\R)$ and each 
$f$ which is holomorphic in a complex neighborhood of the range of $g$, and $\cF$ is \emph{inverse closed} if $1/f \in \cF(U)$ for each 
non-vanishing $f \in \cF(U)$.
That \emph{$\cF$ is derivation closed} means that $f \in \cF(U)$ implies $\p_i f \in \cF(U)$ for all open 
$U \subseteq \R^n$, $n \in \N_{>0}$, and $1 \le i \le n$.
A class $\cF$ of smooth mappings is \emph{quasianalytic} if for each open connected $U \subseteq \R^n$ and each $x \in U$ the 
Borel mapping $\cF(U) \ni f \mapsto (\p^\al f(x))_\al$ is injective.

\section{Weight sequences and \texorpdfstring{$[M]$}{[M]}-ultradifferentiable functions} \label{sec:ws}

\subsection{Weight sequences}
A sequence $M=(M_k) \in \R_{>0}^\N$ 
of positive real numbers is said to 
\begin{basedescript}{\desclabelstyle{\pushlabel}\desclabelwidth{1.2cm}}
    \item[\thetag{M$_{\on{lc}}$}\phantomsection\label{M_lc}] be \emph{log-convex} if $k\mapsto \log M_k$ is convex, 
    i.e., $\forall k : M_k^2 \le M_{k-1} \, M_{k+1}$;
    \item[\thetag{M$_{\on{wlc}}$}\phantomsection\label{M_wlc}] be \emph{weakly log-convex} if $(k!\, M_k)_k$ is log-convex;   
    \item[\thetag{M$_{\on{mg}}$}\phantomsection\label{M_mg}] be of \emph{moderate growth} 
    if $\exists C>0 ~\forall j,k\ge 1 : M_{j+k} \le C^{j+k} M_j \, M_k$;
    \item[\thetag{M$_{\on{dc}}$}\phantomsection\label{M_dc}] be \emph{derivation closed} 
    if $\exists C>0 ~\forall k\ge 1 : M_{k+1} \le C^k M_k$;
    \item[\thetag{M$_{\on{ai}}$}\phantomsection\label{M_ai}] be \emph{almost increasing} if $\exists C>0 ~\forall j \le k : M_j \le C M_k$;
    \item[\thetag{M$_{\on{FdB}}$}\phantomsection\label{M_FdB}] have the \emph{(FdB)-property} if $\exists C>0 ~\forall \al_i \in \N_{>0}, 
    ~\al_1+\cdots+ \al_j=k :  
    M_j M_{\al_1} \cdots M_{\al_j} \le C^{k} M_{k}$;
    \item[\thetag{M$_{\on{qa}}$}\phantomsection\label{M_qa}] be \emph{quasianalytic} 
    if $\sum_{k=1}^\infty (k! M_k)^{-\frac{1}{k}} = \infty$.  
\end{basedescript}
Obviously \thetag{\hyperref[M_lc]{M$_{\on{lc}}$}} implies \thetag{\hyperref[M_wlc]{M$_{\on{wlc}}$}} and 
\thetag{\hyperref[M_mg]{M$_{\on{mg}}$}} implies \thetag{\hyperref[M_dc]{M$_{\on{dc}}$}}. 
If $M$ is log-convex, we further have $M_j M_k \le M_0 M_{j+k}$ for all $j,k$ and $(M_k/M_0)^{\frac{1}{k}}$ is increasing. 
Moreover:

\begin{lemma} \label{FdB}
  For $M \in \R_{>0}^\N$ having the (FdB)-property, each of the following conditions is sufficient:
  \begin{enumerate}[$(1)$]
    \item $M$ is log-convex. 
    \item $M$ is derivation closed and $(M_k^{\frac{1}{k}})_k$ is almost increasing. 
    \item $M_j M_k \le M_1 M_{j+k-1}$ for all $j,k\ge1$.  
  \end{enumerate}
\end{lemma}

\begin{demo}{Proof}
  \thetag{1}
  We show \thetag{\hyperref[M_FdB]{M$_{\on{FdB}}$}} with $C:=\max \{M_1,1\}$ by induction on $k$.
  The assertion is trivial for $k=j$. 
  Assume that $j < k$. Then $\al_i':=\al_i-1 \ge 1$ for some $i$, and
  we have
  \begin{align*}
  M_j M_{\al_1} \cdots M_{\al_j} 
  = M_j M_{\al_1} \cdots M_{\al_i'} \cdots M_{\al_j} \frac{M_{\al_i}}{M_{\al_i'}}
  \le C^{k-1} M_{k-1} \frac{M_k}{M_{k-1}} \le C^{k} M_k,
  \end{align*}
  by induction hypothesis and by \thetag{\hyperref[M_lc]{M$_{\on{lc}}$}}.
   
  \thetag{2} 
  This is proved in more general terms in \ref{thm:{comp}}[\thetag{3} $\Rightarrow$ \thetag{4}] and \ref{thm:fMB}[\thetag{3} $\Rightarrow$ \thetag{4}].
  
  \thetag{3} This is readily seen by iteration.
\qed\end{demo}

For $M,N \in \R_{>0}^\N$ we define: 
\begin{align*}
 M \preceq N \quad :&\Leftrightarrow \quad  \exists C,\rh > 0 ~\forall k : M_k \le C \rh^k N_k \quad \Leftrightarrow \quad 
 \sup_k \Big(\frac{M_k}{N_k}\Big)^{\frac{1}{k}}< \infty  \\
 M \approx N \quad :&\Leftrightarrow \quad M \preceq N \text{ and } N \preceq M \\
 M \lhd N \quad :&\Leftrightarrow \quad \forall \rh>0 ~\exists C>0 ~\forall k : M_k \le C \rh^k N_k 
 \quad \Leftrightarrow \quad \lim_{k \to \infty} \Big(\frac{M_k}{N_k}\Big)^{\frac{1}{k}}=0  
\end{align*}

The following lemma is a variant of \cite[Lemma~6]{Komatsu79b}.

\begin{lemma} \label{Komatsu}
  Let $L,M \in \R_{>0}^\N$ satisfy $L \lhd M$ and $M_k^{\frac{1}{k}} \to \infty$. 
  Then there exist sequences $N^i \in \R_{>0}^\N$, $i=1,2$, satisfying $(N^i_k)^{\frac{1}{k}} \to \infty$ 
    such that $L \le N^1 \lhd N^2 \lhd M$.
\end{lemma}

\begin{demo}{Proof}
  It suffices to show that
  there exists $N^1 \in \R_{>0}^\N$ with $L \le N^1 \lhd M$ and $(N^1_k)^{\frac{1}{k}} \to \infty$; 
  for $N^2=(N^2_k)$ we may then choose $N^2_k := \sqrt{N^1_k M_k}$.

  The sequence $N^1=(N^1_k)$ defined by $N^1_k := \max\{\sqrt{M_k},L_k\}$ is as required: We have
  $L \le N^1 \lhd M$, since
  \[
  \Big(\frac{N^1_k}{M_k}\Big)^{\frac{1}{k}} = \max\Big\{M_k^{-\frac{1}{2k}},\Big(\frac{L_k}{M_k}\Big)^{\frac{1}{k}}\Big\} \to 0
  \] 
  as $M_k^{\frac{1}{k}} \to \infty$ and $L \lhd M$. Moreover, $N^1_k \ge \sqrt{M_k}$
  implies $(N^1_k)^\frac{1}{k} \to \infty$.
\qed\end{demo}

\begin{remark} \label{rmk:Komatsu}
  The lemma remains true if we replace $M_k^{\frac{1}{k}} \to \infty$ by $(k! M_k)^{\frac{1}{k}} \to \infty$ and 
  $(N^i_k)^{\frac{1}{k}} \to \infty$ by
  $(k! N^i_k)^{\frac{1}{k}} \to \infty$; set $N^1_k := \max\{\sqrt{M_k/k!},L_k\}$ in the above proof. 
  But in this case it is unclear if $\varliminf M_k^{\frac{1}{k}}>0$ implies $\varliminf (N^i_k)^{\frac{1}{k}}>0$ which we need in 
  Theorem~\ref{thm:reg}.
\end{remark}

\subsection{Regularizations} Cf.\ \cite{Bang46}, \cite{Mandelbrojt52}, or \cite{Koosis98}.
For $M \in \R_{>0}^\N$ with $(k!M_k)^{\frac{1}{k}} \to \infty$ set 
\[
T_M(t) := \sup_{k \in \N} \frac{t^k}{k! M_k}, ~t > 0, \quad \text{ and } \quad 
M^{\flat(c)}_k := \frac{1}{k!} \sup_{t> 0} \frac{t^k}{T_M(t)}. 
\]
Then $T_M = T_{M^{\flat(c)}}$.
The sequence $(k!M^{\flat(c)}_k)_k$ is the largest log-convex minorant of $(k!M_k)_k$; 
in particular, $M$ is weakly log-convex if and only if $M=M^{\flat(c)}$. 
The condition $(k!M_k)^{\frac{1}{k}} \to \infty$ guarantees that $M_k = M^{\flat(c)}_k$ for infinitely many $k$.

We shall also use 
\[
S_M(t) := \max_{k \le t} \frac{t^k}{k! M_k} \quad \text{ and } \quad M^{\flat(o)}_k  := \frac{1}{k!}\sup_{t \ge k} \frac{t^k}{S_M(t)},
\]
and again have $S_M =S_{M^{\flat(o)}}$.

\begin{lemma} \label{monoton}
  Let $M,N \in \R_{>0}^\N$ satisfy $(k!M_k)^{\frac{1}{k}} \to \infty$ and $(k!N_k)^{\frac{1}{k}} \to \infty$.
  Then $M \preceq N$ implies $M^{\flat(c)} \preceq N^{\flat(c)}$ and
  $M \lhd N$ implies $M^{\flat(c)} \lhd N^{\flat(c)}$.
\end{lemma}

\begin{demo}{Proof}
  For $\rh>0$ set $N^\rh= (N^\rh_k) := (\rh^k N_k)$.  
  Easy computations show
  $T_{N^\rh}(t) = T_N(\frac{t}{\rh})$ and thus $(N^\rh)^{\flat(c)} = (N^{\flat(c)})^\rh$. 
  Both assertions follow immediately. 
\qed\end{demo}

\subsection{\texorpdfstring{$[M]$}{[M]}-ultradifferentiable functions}
Let $M \in \R_{>0}^\N$ and let $U \subseteq \R^n$ be open. 
Define
\begin{align*}
  \cE^{(M)}(U) &:= \Big\{f \in C^\infty(U,\R) : \forall K \subseteq U \text{ compact} ~\forall \rh>0 : \|f\|^M_{K,\rh} < \infty \Big\} \\
  \cE^{\{M\}}(U) &:= \Big\{f \in C^\infty(U,\R) : \forall K \subseteq U \text{ compact} ~\exists \rh>0 : \|f\|^M_{K,\rh} < \infty \Big\} \\
  &\|f\|^M_{K,\rh} := \sup\Big\{\frac{\|f^{(k)}(x)\|_{L^k(\R^n,\R)}}{k!\rh^k M_k}:x\in K,k\in\N\Big\} 
\end{align*}
and endow $\cE^{(M)}(U)$ with its natural Fr\'echet space topology and $\cE^{\{M\}}(U)$ with the projective limit topology over $K$ of the 
inductive limit topology over $\rh$; note that it suffices to take countable limits. 
We write $\cE^{[M]}$ for either $\cE^{(M)}$ or $\cE^{\{M\}}$.   
The elements of $\cE^{[M]}(U)$ are called \emph{$[M]$-ultradifferentiable functions}; an $(M)$/$\{M\}$-ultradifferentiable function is 
said to be of Beurling/Roumieu type, respectively.
For compact $K \subseteq U$ with smooth boundary, 
\[
\cE^M_\rh(K) := \{f \in C^\infty(K) : \|f\|^M_{K,\rh} < \infty \}
\]
is a Banach space, and we have 
\begin{align*}
  \cE^{(M)}(U) = \varprojlim_{K \subseteq U} \varprojlim_{m \in \N} \cE^M_{\frac{1}{m}}(K) 
  \quad \text{and} \quad \cE^{\{M\}}(U) = \varprojlim_{K \subseteq U} \varinjlim_{m \in \N} \cE^M_{m}(K);
\end{align*}
we also set
\begin{align*}
  \cE^{(M)}(K) &:= \Big\{f \in C^\infty(K) : \forall \rh>0 : \|f\|^M_{K,\rh} < \infty \Big\} 
  = \varprojlim_{m \in \N} \cE^M_{\frac{1}{m}}(K)  \\
  \cE^{\{M\}}(K) &:= \Big\{f \in C^\infty(K) : \exists \rh>0 : \|f\|^M_{K,\rh} < \infty \Big\} = \varinjlim_{m \in \N} \cE^M_{m}(K).  
\end{align*}
The definitions work as well for mappings $f: U \to \R^m$, and so we shall use also  
$\cE^{[M]}(U,V)$, $\cE^{[M]}(K,V)$, and $\cE^M_\rh(K,V)$, for open subsets $V \subseteq \R^m$. 

By the Denjoy--Carleman theorem, $\cE^{[M]}$ is quasianalytic if and only if $M^{\flat(c)}$ satisfies 
\thetag{\hyperref[M_qa]{M$_{\on{qa}}$}};
this is in turn equivalent to   
\[
\sum_{k=0}^\infty \frac{M^{\flat(c)}_k}{(k+1) M^{\flat(c)}_{k+1}} = \infty \quad \text{ and } \quad 
\int_1^\infty \frac{\log T_M(t)}{t^2} dt = \infty
\]
For contemporary proofs see for instance \cite[1.3.8]{Hoermander83I}, \cite[19.11]{Rudin87}, and \cite[4.2]{Komatsu73}.

\begin{examples}
  For $s \in \R_{\ge 0}$ the sequence
  $G^s = (G^s_{k}) = ((k!)^s)$ is log-convex and has moderate growth; it is quasianalytic if and only if $s=0$. 
  The elements of $\cE^{\{G^0\}}(U)$ are exactly the real analytic functions $C^\om(U)$ and 
  the elements of $\cE^{(G^0)}(U)$ are exactly the restrictions of entire functions $\cH(\C^n)$. 
  The class $\cE^{\{G^s\}}$ coincides with the Gevrey class $\cG^{1+s}$.
\end{examples}

\begin{lemma} \label{lem:char}
  Let $M \in \R_{>0}^\N$ be weakly log-convex. Then there exists a function 
  $f \in \cE^{\{M\}}_{\on{global}}(\R):=\{f \in C^\infty(\R) : \exists \rh>0 : \|f\|^M_{\R,\rh} < \infty\}$ such that 
  $|f^{(k)}(0)| \ge k! M_k$ for all $k$.
\end{lemma}

Such a function is called a \emph{characteristic $\cE^{\{M\}}$-function}.

\begin{demo}{Proof}
  The complex valued function
  \begin{equation} \label{eq:char}
    g(t) := \sum_{k=0}^\infty \frac{k! M_k}{(2 \mu_k)^k} e^{2 i \mu_k t}, \quad \text{ where } \quad \mu_k := \frac{(k+1) M_{k+1}}{M_k},
  \end{equation}
  belongs to $\cE^{\{M\}}_{\on{global}}(\R,\C)$ and satisfies 
  \begin{align} \label{eq:char1}
    g^{(j)}(0) = i^j h_j, \quad \text{ where } \quad h_j \ge j! M_j,   
  \end{align}
  thus  
  \begin{align*}
    |g^{(j)}(0)| \ge  j! M_j,   
  \end{align*}
  for all $j$; see \cite[Thm.~1]{Thilliez08}.
  Setting $f:= \on{Re} g + \on{Im} g$ we obtain a real valued function with the required properties.
\qed\end{demo}

\begin{proposition} \label{prop:union}
  Let $L, M, N \in \R_{>0}^\N$, let $U \subseteq \R^n$ be open, and let $K \subseteq U$ be compact. We have:
  \begin{enumerate}[$(1)$]
    \item $M \preceq N \Rightarrow \cE^{[M]} \subseteq \cE^{[N]}$ and 
    $M \lhd N \Rightarrow \cE^{\{M\}} \subseteq \cE^{(N)}$ with continuous inclusions. 
    If $M$ is weakly log-convex, then also the converse implications hold; more precisely,
    $\cE^{[M]}(\R) \subseteq \cE^{[N]}(\R) \Rightarrow M \preceq N$ and 
    $\cE^{\{M\}}(\R) \subseteq \cE^{(N)}(\R) \Rightarrow M \lhd N$.
    \item We have
    \[
      \cE^{\{M\}}(U,\R^m) = \bigcap_{M \lhd N} \cE^{(N)}(U,\R^m) = \bigcap_{M \lhd N} \cE^{\{N\}}(U,\R^m).
    \]
    If $M$ is (weakly) log-convex, 
    then the intersections may be taken over all (weakly) log-convex $M \lhd N$.
    \item If $M_k^{\frac{1}{k}} \to \infty$ then 
    \begin{align*}
      \cE^{(M)}(K,\R^m) &= \bigcup_{\substack{L \lhd M \\ L_k^{\frac{1}{k}} \to \infty}} \cE^{\{L\}}(K,\R^m)  
      = \bigcup_{\substack{L \lhd M \\ L_k^{\frac{1}{k}} \to \infty}} \cE^{(L)}(K,\R^m).
    \end{align*}
    If $(k! M_k)^{\frac{1}{k}} \to \infty$ then the unions may be taken over all $L \lhd M$ with $(k! L_k)^{\frac{1}{k}} \to \infty$.
    If $M$ is log-convex and $\tfrac{M_{k+1}}{M_k} \to \infty$ then the unions may be taken over all log-convex $L \lhd M$ with 
    $\tfrac{L_{k+1}}{L_k} \to \infty$.
  \end{enumerate}  
\end{proposition} 

\begin{demo}{Proof}
  \thetag{1} The directions ``$\Rightarrow$'' are clear by definition, see also \cite[2.3]{KMRc}. 
  If $M$ is weakly log-convex, then the implications $\cE^{\{M\}} \subseteq \cE^{\{N\}} \Rightarrow M \preceq N$ and 
  $\cE^{\{M\}} \subseteq \cE^{(N)} \Rightarrow M \lhd N$ follow from the existence
  of a characteristic $\cE^{\{M\}}$-function, see Lemma~\ref{lem:char}.
  That $\cE^{(M)} \subseteq \cE^{(N)}$ implies $M \preceq N$ is shown in \cite[Thm.~2.2]{Bruna80/81} and in more general terms
  in Proposition~\ref{prop:fMincl}. 
  
  \thetag{2} See \cite[2.4 and 8.2]{KMRc}.
  
  \thetag{3} follows from \thetag{1}, Lemma~\ref{Komatsu}, Remark~\ref{rmk:Komatsu}, and \cite[Lemma~6]{Komatsu79b}.
\qed\end{demo}
 
As the elements of $\cE^{\{1\}}(U)$ are exactly the real analytic functions $C^\om(U)$ and 
the elements of $\cE^{(1)}(U)$ are exactly the restrictions of entire functions $\cH(\C^n)$, 
we may conclude:
\begin{enumerate}[$(1)$]
  \item[$(4)$] $C^\om \subseteq \cE^{\{M\}} \Leftrightarrow \cH(\C^n) \subseteq \cE^{(M)}(U) ~\forall U\subseteq \R^n 
  \Leftrightarrow \varliminf M_k^{\frac{1}{k}}>0$
  \item[$(5)$] $C^\om \subseteq \cE^{(M)}$ if and only if $\lim M_k^{\frac{1}{k}} = \infty$.
  \item[$(6)$] $\cE^{[M]}$ is derivation closed if $M$ satisfies \thetag{\hyperref[M_dc]{M$_{\on{dc}}$}}. If $M$ is weakly log-convex, 
  then \thetag{\hyperref[M_dc]{M$_{\on{dc}}$}} is also necessary for $\cE^{[M]}$ being derivation closed; 
  indeed for $M^{+1}=(M^{+1}_k) := (M_{k+1})$ we have $\cE^{[M^{+1}]}(\R) = \{f' : f \in \cE^{[M]}(\R)\}$. 
\end{enumerate}
In particular, if $L \lhd M$ with $\varliminf L_k^{\frac{1}{k}}>0$ then necessarily $\lim M_k^{\frac{1}{k}} = \infty$, by \thetag{1}, 
\thetag{4}, and \thetag{5}. 
 
Note that $\lim \frac{M_{k+1}}{M_k} = \infty$ implies $\lim M_k^{\frac{1}{k}} = \infty$ and thus $C^\om \subseteq \cE^{(M)}$. 
Indeed, there exists $k_0$ with $M_{k_0} \ge 1$, and 
for every $C>0$ there exists $k_1\ge k_0$ so that $M_k \ge C M_{k-1}$ for all $k>k_1$, 
whence 
$M_k^{\frac{1}{k}} \ge M_{k_0}^{\frac{1}{k}} C^{1-\frac{k_0}{k}} \ge C^{\frac{1}{2}}$ as $k >2 k_1$. 
If $M_k^{\frac{1}{k}}$ is increasing, we have also the converse: 
$\lim M_k^{\frac{1}{k}} = \infty$ implies $\lim \frac{M_{k+1}}{M_k} = \infty$.

\begin{lemma}[{\cite[Lemme~3]{Cartan40}}] \label{lem:Cartan}
  Let $M \in \R_{>0}^\N$ and $\la>0$. If $M_0 \le \la^k M_k$ for all $k$ and 
  \[
  |f^{(k)}(t)| \le k! M_k \quad \text{ for all }\quad t \in [-\la,\la], k \in \N,
  \]
  then
  \[
  |f^{(k)}(0)| \le 2 e^k k! M^{\flat(c)}_k \quad \text{ for all }\quad  k \in \N.
  \] 
\end{lemma}

\begin{proposition} \label{prop:reg}
  Let $M \in \R_{>0}^\N$ satisfy $\varliminf M_k^{\frac{1}{k}}>0$ and $M_0=1$, let $K \subseteq \R^n$ be compact, and let 
  $K_\la:=\bigcup_{x \in K} \overline{B_\la(x)}$, $\la>0$, be a $\la$-neighborhood of $K$.
  Then we have $\cE^{\{M\}}(K_\la) \subseteq \cE^{\{M^{\flat(c)}\}}(K)$ via restriction. 
\end{proposition}

\begin{demo}{Proof}
  By the assumption $\varliminf M_k^{\frac{1}{k}}>0$ there exists $\ta>0$ so that $M_k \ge \ta^k$ for all $k$.
  If $f \in \cE^{\{M\}}(K_\la)$, then $C:=\|f\|^M_{K_\la,\rh}<\infty$, where we may assume that $\rh$ is such that $\rh \la \ta \ge 1$. 
  The function $f_{x,v}(t) := f(x + t v)$ satisfies $\|f_{x,v}\|^M_{[-\la,\la],\rh} \le \|f\|^M_{K_\la,\rh}=C$
  for all $x \in K$ and $v \in S^{n-1}$.
  By Lemma~\ref{lem:Cartan}, we have
  \[
    |d_v^k f(x)| = |f_{x,v}^{(k)}(0)| \le 2 C (e \rh)^k k! M^{\flat(c)}_k \quad \text{ for all }\quad x \in K, v \in S^{n-1}, k \in \N,
  \]
  since $(C \rh^k M_k)^{\flat(c)} = C \rh^k M^{\flat(c)}_k$ (see \ref{monoton}). 
  Thus $f|_K \in \cE^{\{M^{\flat(c)}\}}(K)$, see e.g.\ \cite[7.13.1]{KM97}.
\qed\end{demo}

\begin{theorem} \label{thm:reg}
  Let $M \in \R_{>0}^\N$ and let $U \subseteq \R^n$ be open. We have:
  \begin{enumerate}[$(1)$]
    \item If $\varliminf M_k^{\frac{1}{k}}>0$ then $\cE^{\{M\}}(U) = \cE^{\{M^{\flat(c)}\}}(U)$.
    \item If $\lim M_k^{\frac{1}{k}}= \infty$ then $\cE^{(M)}(U) = \cE^{(M^{\flat(c)})}(U)$.
  \end{enumerate}
  Under these assumptions $\cE^{[M]}(U)$ is an algebra.
\end{theorem} 
 
 \thetag{1} is due to \cite[Thm.~I \& Appendix]{Cartan40}.
 
\begin{demo}{Proof}
  \thetag{1} Apply Proposition~\ref{prop:union}\thetag{1} and Proposition~\ref{prop:reg}.
  
  \thetag{2} Proposition~\ref{prop:union}\thetag{1} implies $\cE^{(M^{\flat(c)})}(U) \subseteq \cE^{(M)}(U)$.  
  Conversely, let $K \subseteq U$ be compact and let $K_\la:=\bigcup_{x \in K} \overline{B_\la(x)} \subseteq U$ be a 
  $\la$-neighborhood of $K$ in $U$. By Proposition~\ref{prop:union}\thetag{3}, Proposition~\ref{prop:reg}, 
  and Lemma~\ref{monoton},
  \begin{align*}
  \cE^{(M)}(K_\la) &= \bigcup \cE^{\{L\}}(K_\la)   
  \subseteq \bigcup \cE^{\{L^{\flat(c)}\}}(K)  
  \subseteq \cE^{(M^{\flat(c)})}(K),  
  \end{align*}
  where the unions are taken over all $L \lhd M$ with $L_k^{\frac{1}{k}} \to \infty$.
  As $K$ was arbitrary, we have $\cE^{(M)}(U) \subseteq \cE^{(M^{\flat(c)})}(U)$.

  The supplement is a well-known consequence of weak log-convexity.
\qed\end{demo} 
 
As a consequence $C^\om \subseteq \cE^{\{M\}} = \cE^{(N)}$ is impossible.  
Assume the contrary. Then,
by \ref{prop:union}\thetag{4}\&\thetag{5} and 
Theorem~\ref{thm:reg}, we may assume that $M$ and $N$ are weakly log-convex, and by Proposition~\ref{prop:union}\thetag{1}, 
we have $M \lhd N$. Setting $L=(L_k)$ with $L_k:=\sqrt{M_k N_k}$ we obtain $M \lhd L \lhd N$, and, by Lemma~\ref{monoton}, we may assume 
that $L$ is weakly log-convex. But then $\cE^{\{M\}} \subseteq \cE^{(L)} \subseteq \cE^{\{L\}} \subseteq \cE^{(N)} = \cE^{\{M\}}$ and thus 
$M \approx L \approx N$, a contradiction.

\section{Stability under composition of \texorpdfstring{$\cE^{[M]}$}{E[M]}} \label{sec:wsc}

For $M \in \R_{>0}^\N$ we define $M^\o = (M^\o_k)$ by setting
\[
M^\o_k := \max\{M_jM_{\al_1}\dots M_{\al_j}: 	\al_i\in \N_{>0}, \al_1+\dots+\al_j = k\}, \quad M^\o_0:=1.
\]
Clearly, $M \preceq M^\o$. We have $M^\o \preceq M$ if and only if $M$ has the (FdB)-property.   

\begin{proposition} \label{prop:com} 
  Let $M \in \R_{>0}^\N$ and let $U \subseteq \R^p$, $V \subseteq \R^q$, and $W \subseteq \R^r$ be open.  
  \begin{enumerate}[$(1)$]
    \item If $g \in \cE^{[M]}(U,V)$ and $f \in \cE^{[M]}(V,W)$, then $f\o g \in \cE^{[M^\o]}(U,W)$
    \item If $M$ has the (FdB)-property, then $\cE^{[M]}$ is stable under composition.
  \end{enumerate}
\end{proposition}

\begin{demo}{Proof}
  \thetag{1} Let $K \subseteq U$ be compact. There exist $C_g,\rh_g>0$ (resp.\ for each $\rh_g >0$ there exists $C_g>0$) such that
  \[
  \frac{\|g^{(k)}(x)\|_{L^k(\R^p,\R^q)}}{k!} \le C_g \rh_g^k M_k \quad \text{ for all } x\in K, k \in \N, 
  \]
  and there exist $C_f,\rh_f>0$ (resp.\ for each $\rh_f >0$ there exists $C_f>0$) such that
  \[
  \frac{\|f^{(k)}(y)\|_{L^k(\R^q,\R^r)}}{k!} \le C_f \rh_f^k M_k \quad \text{ for all } y\in g(K), k \in \N.
  \]
  By Fa\`a di Bruno's formula (\cite{FaadiBruno1855} for the 1-dimensional version; 
  the second sum is over all $\al\in \N_{>0}^j$ with $\al_1+\dots+\al_j =k$)
  \begin{align*} 
  &\frac{\|(f\o g)^{(k)}(x)\|_{L^k(\R^p,\R^q)}}{k!} 
  \le \sum_{j\ge 1} \sum_{\al}
  \frac{\|f^{(j)}(g(x))\|_{L^j(\R^q,\R^r)}}{j!}\;\prod_{i=1}^j\;
  \frac{\|g^{(\al_i)}(x)\|_{L^{\al_i}(\R^p,\R^q)}}{\al_i!}
  \\ 
  &\le \sum_{j\ge 1} \sum_{\al} C_f \rh_f^j  C_g^j \rh_g^k M_j\prod_{i=1}^j M_{\al_i}
  \le C_f \rh_g^k \Big(\sum_{j\ge 1} \binom{k-1}{j-1}  (\rh_f  C_g)^j\Big) M^\o_k  \\
  &\le C_f C_g \rh_f (\rh_g (1+\rh_f  C_g))^k M^\o_k.
  \end{align*}
  This implies the assertion in the Roumieu case.
  For the Beurling case,
  let $\ta > 0$ be arbitrary, and choose $\si>0$ such that $\ta=\sqrt{\si}+\si$. 
  If we set $\rh_g = \sqrt{\si}$ and $\rh_f = \sqrt{\si}/C_g$, then $\|f \o g\|^{M^\o}_{K,\ta} < \infty$.
  
  \thetag{2} follows immediately from \thetag{1} and Proposition~\ref{prop:union}\thetag{1}.
\qed\end{demo}
 
We get a nice characterization of stability under composition if we assume that $\cE^{[M]}$ is stable under derivation.

\begin{theorem} \label{thm:comp_dc}
   Let $M \in \R_{>0}^\N$ and assume that $\cE^{[M]}$ is stable under derivation. 
  Consider the following conditions:
  \begin{enumerate}[$(1)$]
    \item $\cE^{[M]}$ is stable under composition.
    \item $\cE^{[M]}$ is holomorphically closed.
    \item $\cE^{[M]}$ is inverse closed.
    \item $(M^{\flat(c)}_k)^{\frac{1}{k}}$ is almost increasing.
    \item $(M^{\flat(o)}_k)^{\frac{1}{k}}$ is almost increasing.
    \item $M^{\flat(c)}$ has the (FdB)-property.
    \item $M^{\flat(o)}$ has the (FdB)-property.
  \end{enumerate}
  If $\varliminf M_k^{\frac{1}{k}}>0$ then all conditions are equivalent in the Roumieu case $\cE^{[M]}=\cE^{\{M\}}$.
  If $\lim M_k^{\frac{1}{k}} = \infty$ then all conditions are equivalent in any case.
\end{theorem}

\begin{demo}{Proof}
Under the assumption $\varliminf M_k^{\frac{1}{k}}>0$ we have $\cE^{\{M\}} = \cE^{\{M^{\flat(c)}\}}$, by Theorem~\ref{thm:reg}. 
The equivalences $\thetag{4} \Leftrightarrow \thetag{5}$ and and $\thetag{6} \Leftrightarrow \thetag{7}$ follow 
from the fact that $\cE^{\{M\}}(I) = \cE^{\{M^{\flat(o)}\}}(I)$ for open intervals $I$, 
see \cite[6.5.1]{Mandelbrojt52}, which implies $M^{\flat(c)} \approx M^{\flat(o)}$, by \cite[Lemma~II]{Siddiqi90}.
Lemma~\ref{FdB} and \ref{prop:union}\thetag{6} imply $\thetag{4} \Rightarrow \thetag{6}$.  

Let us prove the remaining implications in the Roumieu case $\cE^{[M]}=\cE^{\{M\}}$: 
Since $C^\om \subseteq \cE^{\{M\}}$ by \ref{prop:union}\thetag{4}, 
we clearly have $\thetag{1} \Rightarrow \thetag{2} \Rightarrow \thetag{3}$. The 
implication $\thetag{3} \Rightarrow \thetag{5}$ follows from \cite{Siddiqi90}, and
$\thetag{6} \Rightarrow \thetag{1}$ follows from Proposition~\ref{prop:com}. 
Note that $\thetag{3} \Rightarrow \thetag{4}$ is shown in greater generality in the proof of Theorem \ref{thm:{comp}} below.

Now let us assume the stronger condition $\lim M_k^{\frac{1}{k}} = \infty$ and show the  
remaining implications in the Beurling case $\cE^{[M]}=\cE^{(M)}$:
Since $C^\om \subseteq \cE^{(M)}$ by \ref{prop:union}\thetag{5}, we have 
$\thetag{1} \Rightarrow \thetag{2} \Rightarrow \thetag{3}$.
The implication $\thetag{3} \Rightarrow \thetag{4}$ follows from \cite{Bruna80/81} since 
$\cE^{(M)}(\R) = \cE^{(M^{\flat(c)})}(\R)$ is a Fr\'echet algebra, by Theorem~\ref{thm:reg}, and
$\thetag{6} \Rightarrow \thetag{1}$ follows from Proposition~\ref{prop:com}.
\qed\end{demo}

\subsection{Log-convexity is not necessary for stability under composition} 
There exist classes $\cE^{[M]}$ (containing $C^\om$) which are closed under composition 
and there is no log-convex $N \in \R_{>0}^\N$ such that 
$\cE^{[M]} = \cE^{[N]}$. We need the following lemma.

\begin{lemma}\label{Cartanprop3}
Let $M \in \R_{>0}^\N$ be such that $C^\om \subseteq \cE^{[M]}$ (i.e., $\varliminf M_k^{\frac{1}{k}}>0$ in the Roumieu case and 
$\lim M_k^{\frac{1}{k}}=\infty$ in the Beurling case). 
If there exists a log-convex $N \in \R_{>0}^\N$ such that $\cE^{[M]} = \cE^{[N]}$, then the sequence $k_{i+1}/k_i$ is bounded, 
where the $k_i$ are precisely those $k$ with $M_k=M^{\flat(c)}_k$. 
\end{lemma}

\begin{demo}{Proof}
This is a special case of \cite[Appendix Prop.~3]{Cartan40}. For the reader's convenience we give a short proof.
By Theorem~\ref{thm:reg}, we have $\cE^{[M^{\flat(c)}]} = \cE^{[N]}$ and thus $M^{\flat(c)} \approx N$, 
by Proposition~\ref{prop:union}\thetag{1}.  
Since $N$ is weakly log-convex, we have $N \le M^{\flat(c)} \le M$.
Set 
\[
L := 
\begin{pmatrix}
  N_k & k=k_i\\
  +\infty & \text{otherwise}
\end{pmatrix}^{\flat(c)}.
\]
For $M \in \R_{>0}^\N$ consider the graph $\Ga_M := \{(k,\log(k! M_k)): k \in \N\}$.
Then $\Ga_{M^{\flat(c)}}$ and 
$\Ga_L$ lie on piecewise linear curves with vertices $\{(k_i,\log(k_i! M_{k_i})): i \in \N\}$ and
$\{(k_i,\log(k_i! N_{k_i})): i \in \N\}$, respectively. Since $N$ is weakly log-convex and since $\Ga_L$ lies below $\Ga_{M^{\flat(c)}}$, 
we have $N \le L \le M^{\flat(c)} \approx N$ and hence $L \approx N$. 
As $N$ is log-convex, we have,
for $k_i \le k \le k_{i+1}$,
\begin{align*}
\log(k! L_k) &= \tfrac{k_{i+1}-k}{k_{i+1}-k_i} \log(k_i! N_{k_i}) + \tfrac{k-k_i}{k_{i+1}-k_i}\log(k_{i+1}! N_{k_{i+1}})  \\
             &\ge \tfrac{k_{i+1}-k}{k_{i+1}-k_i} \log k_i! + \tfrac{k-k_i}{k_{i+1}-k_i}\log k_{i+1}! + \log N_k,    
\end{align*}
and therefore
\begin{align} \label{eq:Cartan}
\log \Big(\frac{k! L_k}{k! N_k}\Big)^{\frac{1}{k}} &\ge \tfrac{1}{k} \tfrac{k_{i+1}-k}{k_{i+1}-k_i} \log k_i! 
+ \tfrac{1}{k} \tfrac{k-k_i}{k_{i+1}-k_i}\log k_{i+1}! - \tfrac{1}{k} \log k!.
\end{align}  
By Stirling's formula, for $k_{i+1}/k_i=:a_i$ and $k := 2 k_i$ the right-hand side of \eqref{eq:Cartan} is greater than  
\begin{align*} 
\tfrac{1}{2} \tfrac{a_i-2}{a_i-1} (\log k_i - 1) 
+ \tfrac{1}{2} \tfrac{a_i}{a_i-1} (\log a_i +\log k_i - 1)  - \log (2 k_i) 
= \tfrac{1}{2} \tfrac{a_i}{a_i-1} \log a_i - \log 2 -1,
\end{align*}
and so $L \approx N$ implies that $a_i$ is bounded.
\qed\end{demo}

\begin{example} \label{ex:lc}
  Choose $r \in \R_{\ge 4}$. Set $k_i := k_{i-1} \lceil \log (i+1) \rceil$, $i \ge 2$, $k_1:=3$, where $\lceil x \rceil$ denotes the 
  smallest integer $n\ge x$, and define 
  \begin{align*}
    \mu_k = \mu(r)_k :=
    \begin{cases}
       1 & k=1,2\\
       r^{k} & k = k_i  \\
       r^{k_i-1} & k_i < k < k_{i+1}
    \end{cases},     
    \qquad M_k = M(r)_k &:= \frac{1}{k!} \prod_{j=1}^k \mu_j.
  \end{align*}
  Then $M=(M_k)$ is derivation closed, since $\frac{\mu_k}{k} \le r^k$ for all $k$, and   
  $M$ is not weakly log-convex, since $\mu =(\mu_k)$ is not increasing.   
  By construction we have $M_j M_k \le M_1 M_{j+k-1}$ for all $j,k\ge 1$,
  i.e.,
  \[
  \frac{\mu_1}{1}\cdots\frac{\mu_k}{k} \le \frac{\mu_1}{1} \frac{\mu_{j+1}}{j+1}\cdots\frac{\mu_{j+k-1}}{j+k-1}, \quad j,k\ge 1.
  \]
  Indeed, since $\frac{\mu_k}{k}$ is decreasing for $k_i \le k < k_{i+1}$ and 
  since $\frac{\mu_{k_i+1}}{k_i+1} \le \frac{\mu_{k_{i+2}-1}}{k_{i+2}-1}$ for all $i$,
  it suffices to check that, for all $i$,
  \[
  \frac{\mu_{k_{i+1}-1}}{k_{i+1}-1} \frac{\mu_{k_{i+1}}}{k_{i+1}}  \le \frac{\mu_{k_{i+2}-2}}{k_{i+2}-2} \frac{\mu_{k_{i+2}-1}}{k_{i+2}-1} 
  \] 
  which is a straightforward computation.
  By Lemma~\ref{FdB}\thetag{3} and Proposition~\ref{prop:com}, $\cE^{[M]}$ is stable under composition.

  Consider the graph $\Ga_M := \{P_k:=(k,\log(k! M_k)) : k \in \N\}$. The subset $\{P_k : k_i \le k < k_{i+1}\}$ lies 
  on an affine line with slope $(k_i-1) \log r$. The line that connects the two points $P_{k_i-1}$ and $P_{k_i}$ has slope 
  $k_i \log r$, and the line that connects the two points $P_{k_i-1}$ and $P_{k_{i+1}-1}$ has slope 
  $(k_i-1 + (k_{i+1}-k_i)^{-1}) \log r$.
  All these slopes are strictly increasing to infinity in $i$.   
  We may conclude that the graph $\Ga_{M^{\flat(c)}} := \{(k,\log(k! M^{\flat(c)}_k)) : k \in \N\}$ 
  lies on the piecewise linear curve with vertices 
  $\{P_{k_i-1} : i \in \N\}$ and that $\{k_i-1\}$ is precisely the set of $k$ with $M_k = M^{\flat(c)}_k$.  
  
  As $\frac{M_k}{M_{k-1}} = \frac{\mu_k}{k} \to \infty$ we have $M_k^{\frac{1}{k}} \to \infty$ (see the remarks after \ref{prop:union}), and, 
  by Lemma~\ref{Cartanprop3}, there is no log-convex $N \in \R^\N_{>0}$ such that $\cE^{[M]} = \cE^{[N]}$. It is easy to see 
  that the mapping $r \mapsto  \cE^{[M(r)]}$ is injective.
\end{example}

\section{More general spaces of ultradifferentiable functions} \label{sec:wm}

\subsection{Weight matrices} 
A \emph{weight matrix} $\fM = \{M^\la \in \R_{>0}^{\N} : \la \in \La\}$ is a family of 
weakly log-convex sequences 
$M^\la=(M^\la_k)$ satisfying $M^\la_0=1$,
$\lim_{k} (k! M^\la_k)^\frac{1}{k}=\infty$, and
$M^\la\le M^\mu$ if $\la \le \mu$, where $\La$ is a directed partially ordered set.
Let $\sM = \sM(\La)$ be the set of all weight matrices $\fM$ parameterized by the same set $\La$. 
Consider the following conditions:
\begin{basedescript}{\desclabelstyle{\pushlabel}\desclabelwidth{1.5cm}}
  \item[\thetag{$\fM_{\cH}$}\phantomsection\label{fM_H}] 
  $\forall \la \in \La : \varliminf (M^\la_k)^{\frac{1}{k}}>0$.
  \item[\thetag{$\fM_{(C^\om)}$}\phantomsection\label{fM_(Com)}] $\forall \la \in \La : \lim (M^\la_k)^{\frac{1}{k}}=\infty$.
  \item[\thetag{$\fM_{\{C^\om\}}$}\phantomsection\label{fM_{Com}}] $\exists \la \in \La : \varliminf (M^\la_k)^{\frac{1}{k}}>0$.
  \item[\thetag{$\fM_{(\on{dc})}$}\phantomsection\label{fM_(dc)}] $\forall \la \in \La ~\exists \mu \in \La ~\exists C>0 ~\forall k \in \N: 
  M^\mu_{k+1}\le C^{k} M^{\la}_k$.
  \item[\thetag{$\fM_{\{\on{dc}\}}$}\phantomsection\label{fM_{dc}}] $\forall \la \in \La ~\exists \mu \in \La ~\exists C>0 
  ~\forall k \in \N: M^\la_{k+1}\le C^{k} M^{\mu}_k$.
  \item[\thetag{$\fM_{(\on{mg})}$}\phantomsection\label{fM_(mg)}] $\forall \la \in \La ~\exists \mu \in \La ~\exists C>0  
  ~\forall j,k \in \N: 
  M^\mu_{j+k}\le C^{j+k} M^{\la}_j M^{\la}_k$.
  \item[\thetag{$\fM_{\{\on{mg}\}}$}\phantomsection\label{fM_{mg}}] $\forall \la \in \La ~\exists \mu \in \La ~\exists C>0  
  ~\forall j,k \in \N: M^\la_{j+k}\le C^{j+k} M^{\mu}_j M^{\mu}_k$.
  \item[\thetag{$\fM_{(\on{alg})}$}\phantomsection\label{fM_(alg)}] $\forall \la \in \La ~\exists \mu \in \La  ~\exists C>0 
  ~\forall j,k \in \N: M^{\mu}_j M^{\mu}_k \le C^{j+k} M^\la_{j+k}$.
  \item[\thetag{$\fM_{\{\on{alg}\}}$}\phantomsection\label{fM_{alg}}] $\forall \la \in \La ~\exists \mu \in \La  ~\exists C>0 
  ~\forall j,k \in \N: M^{\la}_j M^{\la}_k \le C^{j+k} M^\mu_{j+k}$.
  \item[\thetag{$\fM_{(\on{FdB})}$}\phantomsection\label{fM_(FdB)}] $\forall \la \in \La ~\exists \mu \in \La : (M^\mu)^\o \preceq M^\la$.
  \item[\thetag{$\fM_{\{\on{FdB}\}}$}\phantomsection\label{fM_{FdB}}] $\forall \la \in \La ~\exists \mu \in \La : (M^\la)^\o \preceq M^\mu$.
  \item[\thetag{$\fM_{(\on{L})}$}\phantomsection\label{fM_(L)}] $\forall \la \in \La ~\forall \rh >0 ~\exists \mu \in \La  ~\exists C >0 
  ~\forall k \in \N : \rh^k M^\mu_k \le C M^{\la}_k$.
  \item[\thetag{$\fM_{\{\on{L}\}}$}\phantomsection\label{fM_{L}}] $\forall \la \in \La ~\forall \rh >0 ~\exists \mu \in \La  ~\exists C >0 
  ~\forall k \in \N : \rh^k M^\la_k \le C M^{\mu}_k$.
  \item[\thetag{$\fM_{(\on{BR})}$}\phantomsection\label{fM_(BR)}] $\forall \la \in \La ~\exists \mu \in \La : M^\mu \lhd M^\la$.
  \item[\thetag{$\fM_{\{\on{BR}\}}$}\phantomsection\label{fM_{BR}}] $\forall \la \in \La ~\exists \mu \in \La : M^\la \lhd M^\mu$.
\end{basedescript}
Obviously, \thetag{\hyperref[fM_(Com)]{$\fM_{(C^\om)}$}} $\Rightarrow$ \thetag{\hyperref[fM_H]{$\fM_{\cH}$}} $\Rightarrow$ 
\thetag{\hyperref[fM_{Com}]{$\fM_{\{C^\om\}}$}} 
and \thetag{\hyperref[fM_(mg)]{$\fM_{[\on{mg}]}$}} $\Rightarrow$ \thetag{\hyperref[fM_(dc)]{$\fM_{[\on{dc}]}$}}.
Both conditions \thetag{\hyperref[fM_(alg)]{$\fM_{(\on{alg})}$}} and \thetag{\hyperref[fM_{alg}]{$\fM_{\{\on{alg}\}}$}} are trivially 
satisfied since all $M^\la$ are weakly log-convex, but see Remarks~\ref{rem:fM}.

Henceforth we assume that $\La$ is $\R$ or any ordered subset of $\R$. 
This will enable us to assume that the limits over $\la \in \La$ in the definition of $[\fM]$-ultradifferentiable functions in \ref{ssec:fM}
are countable. Then $\fM$ is in fact an infinite matrix, and the name \emph{weight matrix} is justified. 
On the other hand it is convenient to admit uncountable index sets $\La$. 

\subsection{\texorpdfstring{$[\fM]$}{[fM]}-ultradifferentiable functions} \label{ssec:fM}
Let $\fM$ be a weight matrix, let $U \subseteq \R^n$ be open, and let $K \subseteq U$ be compact. 
We define
\begin{gather*}
\cE^{(\fM)}(K) := \bigcap_{\la \in \La} \cE^{(M^{\la})}(K) \quad 
\text{and} \quad \cE^{\{\fM\}}(K) := \bigcup_{\la \in \La} \cE^{\{M^{\la}\}}(K), \\
\cE^{(\fM)}(U) := \bigcap_{\la \in \La} \cE^{(M^{\la})}(U) \quad 
\text{and} \quad \cE^{\{\fM\}}(U) := \bigcap_{K \subseteq U}\bigcup_{\la \in \La} \cE^{\{M^{\la}\}}(K),
\end{gather*}
and endow these spaces with their natural topologies:
\begin{align*}
  \cE^{(\fM)}(U) := \varprojlim_{\la\in \La} \cE^{(M^{\la})}(U) \quad 
  \text{and} \quad 
  \cE^{\{\fM\}}(U) := \varprojlim_{K \subseteq U} \varinjlim_{\la \in \La} \cE^{\{M^{\la}\}}(K).  
\end{align*}
It is no loss of generality to assume that the limits are countable.
We write $\cE^{[\fM]}$ for either $\cE^{(\fM)}$ or $\cE^{\{\fM\}}$. 
The elements of $\cE^{[\fM]}(U)$ are called \emph{$[\fM]$-ultradifferentiable functions}.
Note that $\cE^{[\fM]}(U)$ forms an algebra, since all $M^\la$ are weakly log-convex. 

We shall use also  
$\cE^{[\fM]}(U,V)$ and $\cE^{[\fM]}(K,V)$, for open subsets $V \subseteq \R^m$.

The inductive limit
\[
\cE^{\{\fM\}}(K,\R^m) = \varinjlim_{\la \in \La} \varinjlim_{\rh>0} \cE^{M^\la}_\rh(K,\R^m) 
=  \varinjlim_{(\la,\rh)} \cE^{M^\la}_\rh(K,\R^m),
\]
where $(\la,\rh) \le (\mu,\si)$ if and only if $\la\le \mu$ and $\rh\le \si$, is a Silva space. 
Indeed, if $\la \le \mu$ and $\rh < \si$ then 
the inclusion 
\[
\cE^{M^\la}_\rh(K,\R^m) \longrightarrow \cE^{M^\mu}_\rh(K,\R^m) \longrightarrow \cE^{M^\mu}_\si(K,\R^m)
\]  
is compact, since the first inclusion is bounded and the second inclusion is compact, by \cite[Prop.~2.2]{Komatsu73}. 

If $\fM$ satisfies \thetag{\hyperref[fM_(L)]{$\fM_{(\on{L})}$}}, respectively \thetag{\hyperref[fM_{L}]{$\fM_{\{\on{L}\}}$}}, we have 
\begin{align} \label{eq:[L]}
  \begin{split}
    \cE^{(\fM)}(K,\R^m)   
  &=  \varprojlim_{(\la,\rh)} \cE^{M^\la}_\rh(K,\R^m) = \varprojlim_{\la} \cE^{M^\la}_1(K,\R^m), 
  \quad \text{respectively} \\
  \cE^{\{\fM\}}(K,\R^m)   
  &=  \varinjlim_{(\la,\rh)} \cE^{M^\la}_\rh(K,\R^m) = \varinjlim_{\la} \cE^{M^\la}_1(K,\R^m)
  \end{split}
\end{align}
as locally convex spaces, where the latter is a Silva space. Indeed, for $1 < \rh$ and by 
\thetag{\hyperref[fM_{L}]{$\fM_{\{\on{L}\}}$}} 
the inclusion
\[
\cE^{M^\la}_1(K,\R^m)  
\longrightarrow \cE^{M^\la}_\rh(K,\R^m) 
\longrightarrow \cE^{M^\mu}_1(K,\R^m)
\]
is compact. If \thetag{\hyperref[fM_(L)]{$\fM_{(\on{L})}$}} then for each $\la \in \La$ and each $\rh>0$ we find $\mu \in \La$ such that
$\cE^{M^\mu}_1(K,\R^m) \subseteq \cE^{M^\la}_\rh(K,\R^m)$ with continuous inclusion.

If $\fM$ satisfies \thetag{\hyperref[fM_(BR)]{$\fM_{(\on{BR})}$}}, respectively \thetag{\hyperref[fM_{BR}]{$\fM_{\{\on{BR}\}}$}}, we have
\begin{align} \label{eq:BR}
\begin{split}
  \cE^{(\fM)}(U,\R^m)   
  &=  \varprojlim_{\la\in \La} \cE^{(M^{\la})}(U,\R^m) =  \varprojlim_{\la\in \La} \cE^{\{M^{\la}\}}(U,\R^m), \quad \text{respectively} \\
  \cE^{\{\fM\}}(K,\R^m)   
  &=  \varinjlim_{\la \in \La} \cE^{\{M^{\la}\}}(K,\R^m) = \varinjlim_{\la \in \La} \cE^{(M^{\la})}(K,\R^m)
\end{split}
\end{align}
as locally convex spaces.

Among the spaces $\cE^{[\fM]}$ we recover the spaces $\cE^{[M]}$ defined by weight sequences, 
if $\fM = \{M\}$ consists just of a single $M \in \R_{>0}^{\N}$,
and the spaces $\cE^{[\om]}$ defined by weight functions, see Corollary~\ref{cor:rep} below. 
We shall see in Theorem~\ref{thm:ext} that in general $\cE^{[\fM]}$ is different from $\cE^{[M]}$ and from $\cE^{[\om]}$.  

\begin{remarks} \label{rem:fM}
  \thetag{1}
  One can replace the condition that the $M^\la \in \fM$ are weakly log-convex by the condition       
  \thetag{\hyperref[fM_H]{$\fM_{\cH}$}} 
  (resp.\ \thetag{\hyperref[fM_(Com)]{$\fM_{(C^\om)}$}}), 
  and work with the log-convex minorants $(M^\la)^{\flat(c)}$ without changing the space $\cE^{\{\fM\}}(U)$ 
  (resp.\ $\cE^{(\fM)}(U)$), 
  see Proposition~\ref{prop:reg} and Theorem~\ref{thm:reg}. 
  Alternatively, assuming \thetag{\hyperref[fM_(alg)]{$\fM_{[\on{alg}]}$}} makes $\cE^{[\fM]}(U)$ into an algebra as well. 
  The condition $M^\la \le M^\mu$ if $\la \le \mu$ may be relaxed to $M^\la \preceq M^\mu$. 

  \thetag{2} Assuming that $(M^\la_k/M^\mu_k)^{\frac{1}{k}}$ is (ultimately) monotonic in $k$ for all $\la,\mu$, 
  we have either $M^\la \approx M^\mu$ for all $\la,\mu$ or $M^\la \lhd M^\mu$ for all $\la < \mu$.
  That is either $\cE^{[\fM]} = \cE^{[M^\la]}$ for all $\la$ or we have the representations in \eqref{eq:BR}. 
\end{remarks}

For $\fM, \fN \in \sM$ we define 
\begin{align*}
 \fM (\preceq) \fN \quad :&\Leftrightarrow \quad \forall \la \in \La ~\exists \mu \in \La : M^\mu \preceq N^\la \\
 \fM \{\preceq\} \fN \quad :&\Leftrightarrow \quad \forall \la \in \La ~\exists \mu \in \La : M^\la \preceq N^\mu  \\
 \fM [\approx] \fN \quad :&\Leftrightarrow \quad \fM [\preceq] \fN \text{ and } \fN [\preceq] \fM \\
 \fM (\preceq\} \fN \quad :&\Leftrightarrow \quad  \exists \la \in \La ~\exists \mu \in \La : M^\la \preceq N^\mu \\
 \fM \{\lhd) \fN \quad :&\Leftrightarrow \quad  \forall \la \in \La ~\forall \mu \in \La : M^\la \lhd N^\mu
\end{align*}

\begin{proposition} \label{prop:fMincl}
  For $\fM,\fN \in \sM$ 
  we have:
  \begin{enumerate}[$(1)$]
    \item $\fM [\preceq] \fN \Rightarrow \cE^{[\fM]} \subseteq \cE^{[\fN]}$ and 
    $\cE^{[\fM]}(\R) \subseteq \cE^{[\fN]}(\R) \Rightarrow \fM [\preceq] \fN$.
    \item $\fM \{\lhd) \fN \Rightarrow \cE^{\{\fM\}} \subseteq \cE^{(\fN)}$ and 
    $\cE^{\{\fM\}}(\R) \subseteq \cE^{(\fN)}(\R) \Rightarrow  \fM \{\lhd) \fN$.
    \item $\fM (\preceq\} \fN \Rightarrow \cE^{(\fM)} \subseteq \cE^{\{\fN\}}$ and 
    $\cE^{(\fM)}(\R) \subseteq \cE^{\{\fN\}}(\R) \Rightarrow  \fM (\preceq\} \fN$.
  \end{enumerate} 
  All inclusions are continuous.
\end{proposition}

\begin{demo}{Proof}
\thetag{1} 
  That $\fM [\preceq] \fN$ implies $\cE^{[\fM]} \subseteq \cE^{[\fN]}$
  is clear by definition. If $\cE^{\{\fM\}}(\R) \subseteq \cE^{\{\fN\}}(\R)$ then  
  $\fM \{\preceq\} \fN$ follows from the existence of characteristic $\cE^{\{M^\la\}}$-functions, by Lemma~\ref{lem:char}. 
  If $\cE^{(\fM)}(\R) \subseteq \cE^{(\fN)}(\R)$ then this inclusion is continuous, by the closed graph theorem since convergence in 
  $\cE^{(\fM)}(\R)$ implies pointwise convergence; here we follow \cite[Thm.~2.2]{Bruna80/81}. 
  Thus 
  for each $\la \in \La$, each compact $I \subseteq \R$, and each $\ta>0$ there exist $\mu \in \La$, $J \subseteq \R$ compact, and constants 
  $C,\rh>0$ such that 
  \[
    \|f\|^{N^\la}_{I,\ta} \le C \|f\|^{M^\mu}_{J,\rh} \quad \text{ for } \quad f \in  \cE^{(\fM)}(\R).
  \] 
  In particular, for $f_t(x) = e^{itx}$ and $\ta=1$, we obtain 
  \[
  T_{N^\la}(t) = \sup_{k \in \N} \frac{t^k}{k! N^\la_k} \le C \sup_{k \in \N} \frac{t^k}{k! \rh^k M^\mu_k} = C T_{M^\mu}(\tfrac{t}{\rh}), 
  \] 
  and thus
  \[
    k! N^\la_k = \sup_{t> 0} \frac{t^k}{T_{N^\la}(t)} \ge \sup_{t> 0} \frac{t^k}{C T_{M^\mu}(\tfrac{t}{\rh})} 
    = k! \frac{\rh^k}{C} M^\mu_k, 
  \]
  that is $\fM (\preceq) \fN$.
  
  \thetag{2} That $\fM \{\lhd) \fN$ implies $\cE^{\{\fM\}} \subseteq \cE^{(\fN)}$ is clear by definition. The converse 
  follows from the existence of characteristic $\cE^{\{M^\la\}}$-functions. 

  \thetag{3} That $\fM (\preceq\} \fN$ implies $\cE^{(\fM)} \subseteq \cE^{\{\fN\}}$ is clear by definition. 
  Conversely, if $\cE^{(\fM)}(\R) \subseteq \cE^{\{\fN\}}(\R)$ then the closed graph theorem 
  (cf.\ \cite[5.4.1]{Jarchow81}) implies that this 
  inclusion is continuous. Indeed $\cE^{(\fM)}(\R)$ is a Fr\'echet space, $\cE^{\{\fN\}}(\R)$ is projective 
  limit of Silva spaces, hence webbed, 
  and convergence implies pointwise convergence. 
  This and Grothendieck's factorization theorem (e.g.\ \cite[24.33]{MeiseVogt97}) imply that 
  for each compact $I \subseteq \R$ there exist $\la \in \La$, $\ta>0$, $\mu \in \La$,   
  $J \subseteq \R$ compact, and constants 
  $C,\rh>0$ such that
  \[
    \|f\|^{N^\la}_{I,\ta} \le C \|f\|^{M^\mu}_{J,\rh} \quad \text{ for } \quad f \in  \cE^{(\fM)}(\R).
  \]
  Applying this to $f_t(x) = e^{itx}$ we obtain, similarly as in \thetag{1},
  \[
    M^\mu_k \le C (\tfrac{\ta}{\rh})^k N^\la_k, 
  \]
  that is $\fM (\preceq\} \fN$.
\qed\end{demo}

We may conclude:
\begin{enumerate}[$(1)$]
  \item[$(4)$] $\cH(\C^n) \subseteq \cE^{(\fM)}(U)$, for all open $U \subseteq \R^n$, if and only if \thetag{\hyperref[fM_H]{$\fM_{\cH}$}}.
  \item[$(5)$] $C^\om \subseteq \cE^{[\fM]}$ if and only if \thetag{\hyperref[fM_(Com)]{$\fM_{[C^\om]}$}}.
  \item[$(6)$] $\cE^{[\fM]}$ is derivation closed if and only if \thetag{\hyperref[fM_(dc)]{$\fM_{[\on{dc}]}$}}.
\end{enumerate}
Note that for $L \in \R_{>0}^\N$ we have $L (\preceq\} \fM$ if and only if $L \{\preceq\} \fM$; in particular,
$\cH(\C^n) \subseteq \cE^{\{\fM\}}(U)$ if and only if $C^\om(U) \subseteq \cE^{\{\fM\}}(U)$, 
for all open $U \subseteq \R^n$.
Moreover:

\begin{corollary}
  Let $M \in \R_{>0}^\N$ with $\lim M_k^{\frac{1}{k}} = \infty$.
  Then there is no $N \in \R_{>0}^\N$ such that $\cE^{(M)}(\R) \subsetneq \cE^{[N]}(\R) \subsetneq \cE^{\{M\}}(\R)$.
\end{corollary}

\begin{demo}{Proof}
  This follows from Proposition~\ref{prop:fMincl} and Theorem~\ref{thm:reg}.
\qed\end{demo}

\begin{remark}
  It is easy to see that $\cE^{\{\fM\}}$ is non-quasianalytic if and only if there is some $\la \in \La$ such that $M^\la$ is 
  non-quasianalytic. Likewise if $\cE^{(\fM)}$ is non-quasianalytic then $M^\la$ is non-quasianalytic for all $\la \in \La$. 
  Intersections $\bigcap_M \cE^{[M]}$, where $M$ runs through a \emph{large} family of non-quasianalytic weakly log-convex 
  weight sequences, can be quasianalytic, see \cite{KMRq} and references therein.  
  But we do not know whether $\cE^{(\fM)}$ can be quasianalytic if all $M^\la$ are non-quasianalytic and $\La$ is restricted 
  to a 1-parameter family (as assumed in this paper).
\end{remark}

\begin{theorem} \label{thm:{comp}}
For a weight matrix $\fM$ satisfying \thetag{\hyperref[fM_{dc}]{$\fM_{\{\on{dc}\}}$}} and \thetag{\hyperref[fM_{Com}]{$\fM_{\{C^\om\}}$}} the following are equivalent:
\begin{enumerate}[$(1)$]
  \item $\cE^{\{\fM\}}$ is stable under composition. 
  \item $\cE^{\{\fM\}}$ is holomorphically closed.
  \item For all $\la\in \La$ there are $\mu\in \La$ and $C>0$ so that 
  $(M^\la_j)^\frac{1}{j} \le C  (M^\mu_k)^\frac{1}{k}$ if $j\le k$.
  \item $\fM$ satisfies \thetag{\hyperref[fM_{FdB}]{$\fM_{\{\on{FdB}\}}$}}. 
\end{enumerate}
\end{theorem}

\thetag{\hyperref[fM_{Com}]{$\fM_{\{C^\om\}}$}} is only needed for $\thetag{1} \Rightarrow \thetag{2}$; \thetag{\hyperref[fM_{dc}]{$\fM_{\{\on{dc}\}}$}} 
is only needed for $\thetag{3} \Rightarrow \thetag{4}$.

\begin{demo}{Proof}
  $\thetag{1} \Rightarrow \thetag{2}$ This is obvious, by \ref{prop:fMincl}\thetag{5}.

  $\thetag{2} \Rightarrow \thetag{3}$ We prove that \thetag{3} holds if $\cE^{\{\fM\}}$ is inverse closed and follow the idea of \cite{Siddiqi90}. 
  Let $\la \in \La$ be fixed and let $g$ be the function in $\cE^{\{M^\la\}}(\R,\C)$ defined by \eqref{eq:char} (with $M=(M_k)$ replaced by $M^\la=(M^\la_k)$).    
  Choose $H>0$ such that $H>1+\sup_{t \in \R} |g(t)|$. 
  We have $H- g \in \cE^{\{M^\la\}}(\R,\C)$, and thus $f:=(H- g)^{-1} \in \cE^{\{\fM\}}(\R,\C)$, as $\cE^{\{\fM\}}(\R,\C)$ is inverse closed, by assumption.
  Thus, there exist $\mu \in \La$ and constants $C,\rh>0$ so that
  \begin{equation} \label{eq:estf}
    \|f\|^{M^\mu}_{[-1,1],\rh} < C.
  \end{equation}
  By Fa\'a di Bruno's formula and using \eqref{eq:char1}, for $k\ge1$, 
 \begin{align*}
    \frac{f^{(k)}(0)}{k!} 
    &= \sum_{j\ge 1} \sum_{\substack{\al_1+\cdots+\al_j=k\\a_\ell>0}} \frac{1}{(H-g(0))^{j+1}} \prod_{\ell=1}^j \frac{g^{(\al_\ell)}(0)}{\al_\ell!} \\
    &= i^k \sum_{j\ge 1} \sum_{\substack{\al_1+\cdots+\al_j=k\\a_\ell>0}} \frac{1}{(H-g(0))^{j+1}} \prod_{\ell=1}^j \frac{h_{\al_\ell}}{\al_\ell!}.
  \end{align*}
  By \eqref{eq:estf},
  \begin{align*}
    C \rh^k M^\mu_k \ge \frac{|f^{(k)}(0)|}{k!}  
    &= \sum_{j\ge 1} \sum_{\substack{\al_1+\cdots+\al_j=k\\a_\ell>0}} \frac{1}{(H-g(0))^{j+1}} \prod_{\ell=1}^j \frac{h_{\al_\ell}}{\al_\ell!} \\
    &\ge \sum_{j\ge 1} \sum_{\substack{\al_1+\cdots+\al_j=k\\a_\ell>0}} \frac{1}{(H-g(0))^{j+1}} \prod_{\ell=1}^j M^\la_{\al_\ell} \\
    &\ge  \frac{1}{(H-g(0))^{k+1}} \prod_{\ell=1}^j M^\la_{\al_\ell}.
  \end{align*}
  In particular, for $\al_1 = \cdots = \al_j=p$, $p \in \N_{>0}$, we have 
  \[
    C_1 \rh_1^{pj} M^\mu_{pj} \ge (M^\la_{p})^j
  \] 
  and hence, for all $j$ and $p$,
  \[
    C_2  (M^\mu_{pj})^{\frac 1 {pj}} \ge (M^\la_{p})^{\frac 1 {p}}.
  \]
  For arbitrary $p\le k$ choose $j$ so that $jp \le k < (j+1)p$. Then 
  \[
    (M^\mu_{k})^{\frac 1 {k}} \ge (M^\mu_{jp})^{\frac 1 {jp}} \frac{(jp)!^{\frac 1 {jp}}}{k!^{\frac 1 {k}}} \ge 
    C_2^{-1} (M^\la_{p})^{\frac 1 {p}} \frac{(jp)!^{\frac 1 {jp}}}{k!^{\frac 1 {k}}} \ge C_2^{-1} (M^\la_{p})^{\frac 1 {p}},
  \]
  since $(k! M^\mu_k)^{1/k}$ is non-decreasing.

  $\thetag{3} \Rightarrow \thetag{4}$   
  By \thetag{\hyperref[fM_{dc}]{$\fM_{\{\on{dc}\}}$}}, for $\la \in \La$ 
  there exist $\mu\in \La$ and $D>0$ so that $M^\la_{k+1} \le D^k M^\mu_k$ for all $k\ge 1$. 
  The assumption implies that there is $\nu \in \La$ so that    
  $M^\mu_{\be_1} \cdots M^\mu_{\be_j} \le C^k M^\nu_{k}$ for all 
  $\be_i \in \N_{>0}$ with $\be_1 + \cdots +\be_j=k$.  
  Let $I := \{i : \al_i \ge 2\}$ and set $\al_i':=\al_i-1$.
  Then, as $\mu \ge \la$, 
  \begin{align*}
  M^\la_j M^\la_{\al_1} \cdots M^\la_{\al_j} &= M^\la_j (M^\la_1)^{j-|I|} \prod_{i \in I} M^\la_{\al_i} \le D^{k-j} M^\la_j 
  (M^\la_1)^{j-|I|} \prod_{i \in I} M^\mu_{\al_i'} \\
  &\le D^{k-j} (M^\la_1)^{j-|I|} C^{k} M^\nu_{k} \le \tilde C^{k} M^\nu_k,
  \end{align*}  
  which shows \thetag{4}.

  $\thetag{4} \Rightarrow \thetag{1}$ 
  Let $g \in \cE^{\{\fM\}}(U,V)$ and $f \in \cE^{\{\fM\}}(V,W)$, 
  for open subsets $U \subseteq \R^p$, $V \subseteq \R^q$, $W \subseteq \R^r$, and let $K \subseteq U$ be compact. 
  By definition, there exist $\la_i\in \La$, $i=1,2$, such that 
  $g \in \cE^{\{M^{\la_1}\}}(K,V)$ and $f \in \cE^{\{M^{\la_2}\}}(g(K),W)$, and there exists $\la \ge \la_i$, $i=1,2$. 
  By \thetag{\hyperref[fM_{FdB}]{$\fM_{\{\on{FdB}\}}$}}, there exists $\mu\in \La$ such that 
  $(M^\la)^\o \preceq M^\mu$, and thus, 
  by Proposition~\ref{prop:com}, we have $f \o g \in \cE^{\{M^\mu\}}(K,W)$ which implies the assertion.
\qed\end{demo}

\begin{theorem} \label{thm:fMB}
For a weight matrix $\fM$ satisfying \thetag{\hyperref[fM_(dc)]{$\fM_{(\on{dc})}$}} and \thetag{\hyperref[fM_H]{$\fM_{\cH}$}} the following are equivalent:
\begin{enumerate}[$(1)$]
  \item $\cE^{(\fM)}$ is stable under composition. 
  \item $\cE^{(\fM)}$ is holomorphically closed.
  \item For all $\la\in \La$ there are $\mu\in \La$ and $C>0$ so that 
  $(M^\mu_j)^\frac{1}{j} \le C  (M^\la_k)^\frac{1}{k}$ if $j\le k$.
  \item $\fM$ satisfies \thetag{\hyperref[fM_(FdB)]{$\fM_{(\on{FdB})}$}}. 
\end{enumerate}
\end{theorem}

\thetag{\hyperref[fM_H]{$\fM_{\cH}$}} is only needed for $\thetag{1} \Rightarrow \thetag{2}$; \thetag{\hyperref[fM_(dc)]{$\fM_{(\on{dc})}$}} 
is only needed for $\thetag{3} \Rightarrow \thetag{4}$.

\begin{demo}{Proof}
  $\thetag{1} \Rightarrow \thetag{2}$ This is obvious, by \ref{prop:fMincl}\thetag{4}.
  
  $\thetag{2} \Rightarrow \thetag{3}$ We follow \cite{Bruna80/81}. Since all $M^\la$ are weakly log-convex, $\cE^{(\fM)}(\R)$ is a 
  Fr\'echet algebra which is locally m-convex, by \cite{MitiaginRolewiczZelazko62}, i.e., $\cE^{(\fM)}(\R)$ has an equivalent seminorm
  system $\{p\}$ such that $p(fg) \le p(f)p(g)$ for all $f,g \in \cE^{(\fM)}(\R)$. 
  So for each $\la\in \La$, compact $K \subseteq \R$, and $\rh>0$ there exist $p$, $\mu\in \La$, compact $L \subseteq \R$, $\si>0$ 
  and constants $C,D >0$ such that
  \[ 
  \|f^m\|^{M^{\la}}_{K,\rh} \le C p(f^m) \le C (p(f))^m \le C D^m (\|f\|^{M^\mu}_{L,\si})^m, \quad f \in \cE^{(\fM)}(\R), m \in \N,
  \]
  in particular, for $f_t(x) = e^{itx}$ and $\rh=1$, we find
  \[
  T_{M^\la}(mt) \le C D^m (T_{M^\mu}(\tfrac{t}{\si}))^m.
  \]
  Let $j \le k$ and suppose that $k = jl$ with $l \in \N$. We have, for some constant $\tilde C$, 
  \[
  (T_{M^\la}(t))^\frac{1}{k} = (T_{M^\la}(l\tfrac{t}{l}))^\frac{1}{k}  
  \le C^\frac{1}{k} D^\frac{1}{j} (T_{M^\mu}(\tfrac{t}{\si l}))^\frac{1}{j} 
  \le \tilde C (T_{M^\mu}(\tfrac{t}{\si l}))^\frac{1}{j},
  \]
  thus
  \begin{align*}
    (k!M^\la_k)^\frac{1}{k} =  \sup_{t>0} \frac{t}{(T_{M^\la}(t))^\frac{1}{k}}
    \ge  \sup_{t>0} \frac{t}{\tilde C (T_{M^\mu}(\tfrac{t}{\si l}))^\frac{1}{j}} 
    =  \frac{\si l}{\tilde C} (j!M^\mu_j)^\frac{1}{j}.
  \end{align*}
  In general choose $l \in \N$ such that $lj \le k < (l+1)j$. Then, as $(k!M^\la_k)^\frac{1}{k}$ is increasing,  
  \[
  (k!M^\la_k)^\frac{1}{k} \ge ((lj)!M^\la_{lj})^\frac{1}{lj} \ge \frac{\si l}{\tilde C} (j!M^\mu_j)^\frac{1}{j} \ge 
  \frac{\si (l+1)}{2 \tilde C} (j!M^\mu_j)^\frac{1}{j},
  \]
  and, by Stirling's formula, there is a constant $\bar C>0$ such that $(M^\mu_j)^\frac{1}{j} \le \bar C (M^\la_k)^\frac{1}{k}$ for 
  all $j\le k$.

  $\thetag{3} \Rightarrow \thetag{4}$  
  The assumption implies that     
  $M^\mu_{\be_1} \cdots M^\mu_{\be_j} \le C^k M^\la_{k}$ for all 
  $\be_i \in \N_{>0}$ with $\be_1 + \cdots +\be_j=k$. 
  By \thetag{\hyperref[fM_(dc)]{$\fM_{(\on{dc})}$}}, 
  there exist $\nu\in \La$ and $D>0$ so that $M^\nu_{k+1} \le D^k M^\mu_k$ for all $k\ge 1$. 
  Let $I := \{i : \al_i \ge 2\}$ and set $\al_i':=\al_i-1$.
  Then, as $\nu \le \mu$, 
  \begin{align*}
  M^\nu_j M^\nu_{\al_1} \cdots M^\nu_{\al_j} &= M^\nu_j (M^\nu_1)^{j-|I|} \prod_{i \in I} M^\nu_{\al_i} \le D^{k-j} M^\nu_j 
  (M^\nu_1)^{j-|I|} \prod_{i \in I} M^\mu_{\al_i'} \\
  &\le D^{k-j} (M^\nu_1)^{j-|I|} C^{k} M^\la_{k} \le \tilde C^{k} M^\la_k,
  \end{align*}  
  which shows \thetag{4}.
  
  $\thetag{4} \Rightarrow \thetag{1}$ 
  Let $g \in \cE^{(\fM)}(U,V)$ and $f \in \cE^{(\fM)}(V,W)$, 
  for open subsets $U \subseteq \R^p$, $V \subseteq \R^q$, $W \subseteq \R^r$, and let $K \subseteq U$ be compact. 
  By definition, for each $\mu\in \La$ we have 
  $g \in \cE^{(M^{\mu})}(K,V)$ and $f \in \cE^{(M^{\mu})}(g(K),W)$. 
  By \thetag{\hyperref[fM_(FdB)]{$\fM_{(\on{FdB})}$}} and by Proposition~\ref{prop:com}, 
  we obtain $f \o g \in \cE^{(M^{\la})}(K,W)$ for each $\la \in \La$ which implies the assertion. 
\qed\end{demo}

\subsection{Composition operators}

Let $\fM$ be a weight matrix.
If $\fM$ satisfies \thetag{\hyperref[fM_(FdB)]{$\fM_{[\on{FdB}]}$}}, we may consider the nonlinear composition operators
\begin{gather*}
 \on{comp}^{[\fM]} : \cE^{[\fM]}(\R^p,\R^q) \times \cE^{[\fM]}(\R^q,\R^r) \to \cE^{[\fM]}(\R^p,\R^r) : (g,f) \mapsto f \o g \\
 \cE^{[\fM]}(\R^p,f) : \cE^{[\fM]}(\R^p,\R^q) \to \cE^{[\fM]}(\R^p,\R^r) : g \mapsto f \o g, \quad f \in \cE^{[\fM]}(\R^q,\R^r),
\end{gather*}
by Theorem~\ref{thm:{comp}} and Theorem~\ref{thm:fMB}.

\begin{theorem} \label{thm:compC0}
  We have:
  \begin{enumerate}[$(1)$]
    \item If $\fM$ satisfies \thetag{\hyperref[fM_(FdB)]{$\fM_{(\on{FdB})}$}}, then $\on{comp}^{(\fM)}$ is continuous.
    \item If $\fM$ satisfies \thetag{\hyperref[fM_{FdB}]{$\fM_{\{\on{FdB}\}}$}}, then $\cE^{\{\fM\}}(\R^p,f)$, 
    for $f \in \cE^{\{\fM\}}(\R^q,\R^r)$, is continuous 
    and $\on{comp}^{\{\fM\}}$ is sequentially continuous.
  \end{enumerate}
\end{theorem}

\begin{demo}{Proof} 
  We follow \cite{AppellNazarovZabrejko91}
  and subdivide the proof into several claims.
  
  \begin{claim} \label{cl:1}
    If $\fM$ satisfies \thetag{\hyperref[fM_(FdB)]{$\fM_{[\on{FdB}]}$}}, then $\on{comp}^{[\fM]}$ is bounded.
  \end{claim}
  
  We treat the cases $\cE^{(\fM)}$ and $\cE^{\{\fM\}}$ separately. 
    
  ($\cE^{[\fM]}=\cE^{\{\fM\}}$) 
  Let $\cB_1 \subseteq \cE^{\{\fM\}}(\R^p,\R^q)$ and $\cB_2 \subseteq \cE^{\{\fM\}}(\R^q,\R^r)$ be bounded subsets.
  Let $K \subseteq \R^p$ be an arbitrary, but fixed, compact subset. Then $\cB_1$ is bounded in $\cE^{\{\fM\}}(K,\R^q)$.  
  Since the inductive limit $\cE^{\{\fM\}}(K,\R^q) = \varinjlim_{(\la,\rh)} \cE^{M^\la}_\rh(K,\R^q)$ is regular, see \ref{ssec:fM},
  $\cB_1$ is contained and bounded in some step $\cE^{M^{\la_1}}_{\rh_1}(K,\R^q)$, i.e., there exist $\la_1 \in \La$ and $\rh_1>0$ such that 
  $\sup_{g \in \cB_1} \|g\|^{M^{\la_1}}_{K,{\rh_1}} < \infty$. In particular, the closure 
  \begin{equation} \label{eq:L}
   L:=\overline{\bigcup_{g \in \cB_1} g(K)}  
  \end{equation}
  is a compact subset of $\R^q$, and $\cB_2$ is bounded in $\cE^{\{\fM\}}(L,\R^r) = \varinjlim_{(\la,\rh)} \cE^{M^\la}_\rh(L,\R^r)$.
  So there exist $\la_2 \in \La$ and $\rh_2>0$ such that 
  $\sup_{f \in \cB_2} \|f\|^{M^{\la_2}}_{L,{\rh_2}} < \infty$.
  For $\la := \max \la_i$, we have 
  \begin{equation} \label{eq:C1C2}
  C_1:=\sup_{g \in \cB_1} \|g\|^{M^\la}_{K,\rh_1} < \infty \quad \text{ and } \quad 
  C_2 := \sup_{f \in \cB_2} \|f\|^{M^\la}_{L,\rh_2}  < \infty.  
  \end{equation}
  By the proof of Proposition~\ref{prop:com} we find that 
  \begin{equation} \label{eq:prop31}
    \sup_{(g,f) \in \cB_1 \times \cB_2} \|f \o g\|^{(M^\la)^\o}_{K,\si} 
    \le C_1 C_2 \rh_2 < \infty, \quad \text{ with } \si :=\rh_1(1+\rh_2 C_1),
  \end{equation}  
  and by \thetag{\hyperref[fM_{FdB}]{$\fM_{\{\on{FdB}\}}$}} there exist $\mu \in \La$ and $C>0$ such that 
  \begin{equation} \label{eq:ccl}
    \sup_{(g,f) \in \cB_1 \times \cB_2} \|f \o g\|^{M^\mu}_{K,C \si} 
    \le \sup_{(g,f) \in \cB_1 \times \cB_2} \|f \o g\|^{(M^\la)^\o}_{K,\si} < \infty.
  \end{equation}
  Since $K$ was arbitrary, $\on{comp}^{\{\fM\}}(\cB_1 \times \cB_2)$ is bounded in $\cE^{\{\fM\}}(\R^p,\R^r)$. 
  
  ($\cE^{[\fM]}=\cE^{(\fM)}$) 
  Let $\cB_1 \subseteq \cE^{(\fM)}(\R^p,\R^q)$ and $\cB_2 \subseteq \cE^{(\fM)}(\R^q,\R^r)$ be bounded. 
  Let $\mu \in \La$, let $K \subseteq \R^p$ be compact, and let $\ta>0$.
  By \thetag{\hyperref[fM_(FdB)]{$\fM_{(\on{FdB})}$}}, 
  we find $\la \in \La$ and $C>0$ such that $(M^\la)^\o_k \le C^k M^\mu_k$ for all $k$.
  Choose $\rh>0$ so that $\ta/C = \sqrt{\rh}+\rh$ and set $\rh_1 = \sqrt{\rh}$.
  Let $C_1$ be defined by \eqref{eq:C1C2};
  $\cB_1$ is bounded in $\cE^{M^{\la}}_{\rh_1}(K,\R^q)$. 
  Set $\rh_2 = \sqrt{\rh}/C_1$ and let $C_2$ be defined by \eqref{eq:C1C2};
  $\cB_2$ is bounded in $\cE^{M^{\la}}_{\rh_2}(L,\R^r)$, where $L$ is defined by 
  \eqref{eq:L}. 
  As before we may conclude \eqref{eq:prop31} and \eqref{eq:ccl}, where $\si=\ta/C$,
  which completes the proof of the claim.
  
  \begin{claim} \label{cl:2}
    If $\fM$ satisfies \thetag{\hyperref[fM_(FdB)]{$\fM_{[\on{FdB}]}$}}, then $\on{comp}^{[\fM]}$ is sequentially continuous.
  \end{claim}  
  
  ($\cE^{[\fM]}=\cE^{\{\fM\}}$)
  Let $(g_n,f_n) \to (g,f)$ in $\cE^{\{\fM\}}(\R^p,\R^q) \times \cE^{\{\fM\}}(\R^q,\R^r)$. Then the sets 
  $\cB_1 := \{g_n : n \in \N\}$, $\cB_2:=\{f_n : n \in \N\}$, 
  and $\{f_n \o g_n : n \in \N\} \subseteq \on{comp}^{\{\fM\}}(\cB_1 \times \cB_2)$
  are bounded, by Claim~\ref{cl:1}.
  Let $K \subseteq \R^p$ be an arbitrary, but fixed, compact subset, and let $L$ be given by \eqref{eq:L}.
  By regularity of the inductive limit $\cE^{\{\fM\}}(K,\R^r) = \varinjlim_{(\la,\rh)} \cE^{M^\la}_\rh(K,\R^r)$, the set 
  $\{f_n \o g_n : n \in \N\}$ is contained and bounded in some step $\cE^{M^\la}_{\rh}(K,\R^r)$
  and hence is
  precompact in $\cE^{M^\mu}_{\si}(K,\R^r)$, where  $\la \le \mu$ and $\rh < \si$, see \ref{ssec:fM}, and so it has an accumulation point 
  $h \in \cE^{M^\mu}_{\si}(K,\R^r)$.  
  It is well-known that composition of continuous mappings, i.e., $\on{comp}^{0} : C^0(K,L) \times C^0(L,\R^r) \to C^0(K,\R^r)$, is 
  continuous, see e.g.\ \cite[Thm.~3.4.2]{Engelking89}, 
  and thus $f_n \o g_n \to f \o g$ in $C^0(K,\R^r)$.  
  It follows that $f \o g = h$. As $K$ was arbitrary the assertion follows.
  
  ($\cE^{[\fM]}=\cE^{(\fM)}$)
  The proof is analogous; note that here $\{f_n \o g_n : n \in \N\}$ is precompact in every step $\cE^{M^\mu}_{\si}(K,\R^r)$.
  
  \begin{claim} \label{cl:3}
    If $\fM$ satisfies \thetag{\hyperref[fM_(FdB)]{$\fM_{(\on{FdB})}$}}, then $\on{comp}^{(\fM)}$ is continuous.
  \end{claim} 
  
  This follows from Claim~\ref{cl:2}, since $\cE^{(\fM)}(\R^p,\R^q) \times \cE^{(\fM)}(\R^q,\R^r)$ is metrizable.
  
  \begin{claim} \label{cl:4}
    If $\fM$ satisfies \thetag{\hyperref[fM_{FdB}]{$\fM_{\{\on{FdB}\}}$}}, then $\cE^{\{\fM\}}(\R^p,f)$ is continuous.
  \end{claim}
  
  Arguments similar as in the proof of Claim~\ref{cl:2} show that the restricted mapping 
  $\cE^{\{\fM\}}(K,f) : \cE^{\{\fM\}}(K,\R^q) \to \cE^{\{\fM\}}(K,\R^r)$ is sequentially continuous, thus continuous, 
  for each compact subset $K \subseteq \R^p$, 
  since $\cE^{\{\fM\}}(K,\R^q)$ is sequential, by \ref{ssec:fM} and e.g.\ \cite[8.5.28]{Perez-CarrerasBonet87}. The projective structure of 
  $\cE^{\{\fM\}}(\R^p,\R^q) = \varprojlim_{K} \cE^{\{\fM\}}(K,\R^q)$ implies that the mapping
  $\cE^{\{\fM\}}(\R^p,f) : \cE^{\{\fM\}}(\R^p,\R^q) \to \cE^{\{\fM\}}(\R^p,\R^r)$ is continuous.
\qed\end{demo}

\begin{corollary}
  Let $M \in \R_{>0}^\N$ satisfy \thetag{\hyperref[M_FdB]{M$_{\on{FdB}}$}}.
  Then $\on{comp}^{(M)}$ is continuous, $\cE^{\{M\}}(\R^p,f)$, for $f \in \cE^{\{M\}}(\R^q,\R^r)$, is continuous, and $\on{comp}^{\{M\}}$ is sequentially continuous.  
\end{corollary}

\begin{demo}{Proof}
  This is a special case of Theorem~\ref{thm:compC0}; weak log-convexity of $M$ is not needed here.
\qed\end{demo}

\begin{remark} \label{rem:fc}
  If $M$ additionally has moderate growth, then the mapping 
  $\on{comp}^{[M]}$ is 
  even $\cE^{[M]}$ which is a consequence of the $\cE^{[M]}$-exponential law, see \cite[5.5]{KMRu}. 
  We expect that more generally 
  $\on{comp}^{[\fM]}$ is $\cE^{[\fM]}$, if $\fM$ satisfies \thetag{\hyperref[fM_(FdB)]{$\fM_{[\on{FdB}]}$}} and  
  \thetag{\hyperref[fM_(mg)]{$\fM_{[\on{mg}]}$}}. 
  This is work in progress and will appear in a forthcoming paper. 
\end{remark}

\section{Weight functions and \texorpdfstring{$[\om]$}{[om]}-ultradifferentiable functions} \label{sec:wf}

\subsection{Weight functions} \label{ssec:wf}
Let $\sW$ be the set of all continuous increasing functions $\om: [0,\infty) \to [0,\infty)$ with $\om|_{[0,1]}=0$, 
$\lim_{t \to \infty} \om(t) = \infty$, and so that the following assumptions \thetag{\hyperref[om_1]{$\om_1$}}, 
\thetag{\hyperref[om_2]{$\om_2$}}, and \thetag{\hyperref[om_3]{$\om_3$}}
are satisfied:
\begin{basedescript}{\desclabelstyle{\pushlabel}\desclabelwidth{0.7cm}} 
  \item[\thetag{$\om_1$}\phantomsection\label{om_1}] $\om(2t)=O(\om(t))$ as $t\to \infty$.
  \item[\thetag{$\om_2$}\phantomsection\label{om_2}] $\log(t)=o(\omega(t))$ as $t\to \infty$.
  \item[\thetag{$\om_3$}\phantomsection\label{om_3}] $\vh : t \mapsto \om(e^t)$ is convex on $[0,\infty)$.
\end{basedescript}  
Occasionally, we shall also consider the following conditions:
\begin{basedescript}{\desclabelstyle{\pushlabel}\desclabelwidth{0.7cm}}
  \item[\thetag{$\om_4$}\phantomsection\label{om_4}] $\om(t)=O(t)$ as $t\to \infty$.
  \item[\thetag{$\om_5$}\phantomsection\label{om_5}] $\om(t)=o(t)$ as $t \to \infty$.
  \item[\thetag{$\om_6$}\phantomsection\label{om_6}] $\exists H\ge 1 ~\forall t\ge 0 : 2 \om(t) \le \om(H t)+H$.
  \item[\thetag{$\om_7$}\phantomsection\label{om_7}] $\exists C>0 ~\exists t_0>0 ~\forall \la\ge 1 ~\forall t\ge t_0 
  : \om(\la t) \le C \la \om(t)$.
  \item[\thetag{$\om_8$}\phantomsection\label{om_8}] $\exists C>0 ~\exists H>0 ~\forall t\ge 0 : \om(t^2) \le C \om(H t) +C$. 
\end{basedescript}
Then $\sW$ forms an abelian semigroup with respect to pointwise addition, which also preserves all conditions 
\thetag{\hyperref[om_4]{$\om_4$}}--\thetag{\hyperref[om_8]{$\om_8$}}.

For $\om \in \sW$
the \emph{Young conjugate} of $\vh$, given by 
\[
\vh^*(t):=\sup \{st-\vh(s) : s \ge 0\}, \quad t \ge 0,
\] 
is convex, increasing, and satisfies $\vh^*(0)=0$, $\vh^{**}=\vh$,  
and $\lim_{t\to \infty}\frac{t}{\vh^*(t)}=0$. 
Moreover,
the functions $t \mapsto \tfrac{\vh(t)}{t}$ and $t \mapsto \tfrac{\vh^*(t)}{t}$ are increasing. 
Cf.\ \cite{BMT90}. 
Convexity of $\vh^*$ and $\vh^*(0)=0$ implies
\begin{equation} \label{eq:cx0}
\vh^*(t)+\vh^*(s) \le \vh^*(t+s) \le \tfrac{1}{2}\vh^*(2t)+ \tfrac{1}{2}\vh^*(2s) , \quad t,s\ge 0. 
\end{equation}

Note that $\om(t) := \max\{0,(\log t)^s\}$, $s >1$, belongs to $\sW$ and satisfies all listed conditions except 
\thetag{\hyperref[om_6]{$\om_6$}}. 

For $\om, \si \in \sW$ we define:
\begin{align*}
 \om \preceq \si \quad :&\Leftrightarrow \quad  \si(t)=O(\om(t)) \text{ as } t\to \infty  \\
 \om \approx \si \quad :&\Leftrightarrow \quad \om \preceq \si \text{ and } \si \preceq \om \\
 \om \lhd \si \quad :&\Leftrightarrow \quad  \si(t)=o(\om(t)) \text{ as } t\to \infty   
\end{align*}

\subsection{\texorpdfstring{$[\om]$}{[om]}-ultradifferentiable functions}  \label{subsec:om}
Let $\om \in \sW$ and let $U \subseteq \R^n$ be open. 
Define
\begin{align*}
  \cE^{(\om)}(U) 
  &:= \Big\{f \in C^\infty(U,\R) : 
  \forall K \subseteq U \text{ compact} ~\forall \rh > 0 : \|f\|^\om_{K,\rh} < \infty \Big\} \\
  \cE^{\{\om\}}(U) 
  &:= \Big\{f \in C^\infty(U,\R) : \forall K \subseteq U \text{ compact} ~\exists \rh>0 : \|f\|^\om_{K,\rh} < \infty \Big\} \\
  &\|f\|^\om_{K,\rh} := \sup\Big\{\|f^{(k)}(x)\|_{L^k(\R^n,\R)} \exp(- \tfrac{1}{\rh} \vh^*(\rh k)) : x\in K,k\in\N\Big\} 
\end{align*}
and endow $\cE^{(\om)}(U)$ with its natural Fr\'echet space topology and $\cE^{\{\om\}}(U)$ with 
the projective limit topology over $K$ of the 
inductive limit topology over $\rh$; note that it suffices to take countable limits. 
We write $\cE^{[\om]}$ for either $\cE^{(\om)}$ or $\cE^{\{\om\}}$. 
The elements of $\cE^{[\om]}(U)$ are called \emph{$[\om]$-ultradifferentiable functions}; 
an $(\om)$/$\{\om\}$-ultradifferentiable function is said to be of Beurling/Roumieu type, respectively. For compact $K \subseteq U$ with 
smooth boundary, set
\begin{align*}
  \cE^{\om}_\rh(K) &:= \Big\{f \in C^\infty(K) : \|f\|^\om_{K,\rh} < \infty \Big\} \\
  \cE^{(\om)}(K) &:= \Big\{f \in C^\infty(K) : \forall \rh > 0 : \|f\|^\om_{K,\rh} < \infty \Big\} 
  = \varprojlim_{m \in \N} \cE^{\om}_{\frac{1}{m}}(K) \\
  \cE^{\{\om\}}(K) &:= \Big\{f \in C^\infty(K) : \exists \rh > 0 : \|f\|^\om_{K,\rh} < \infty \Big\}
  = \varinjlim_{m \in \N} \cE^{\om}_{m}(K).  
\end{align*}
We shall also use   
$\cE^{[\om]}(U,V)$, $\cE^{[\om]}(K,V)$, and $\cE^{\om}_\rh(K,V)$ for open subsets $V \subseteq \R^m$.

Note that 
$\cE^{[\om]}$ is quasianalytic if and only if 
\[
\int_1^\infty \frac{\om(t)}{t^2} dt = \infty
\]
(e.g., by Corollary~\ref{cor:omqa} and Theorem~\ref{thm:rep} below),
and in this case we say that $\om$ is quasianalytic.

\begin{examples}
  For $s \in \R_{\ge 0}$ the weight function
  $\ga^s(t) = t^{\tfrac{1}{1+s}}$ 
  has all properties listed in \ref{ssec:wf} except \thetag{\hyperref[om_8]{$\om_8$}} and \thetag{\hyperref[om_5]{$\om_5$}} if $s=0$;
  it is quasianalytic if and only if $s=0$. 
  The elements of $\cE^{\{\ga^0\}}(U)$ are exactly the real analytic functions $C^\om(U)$ and 
  the elements of $\cE^{(\ga^0)}(U)$ are exactly the restrictions of entire functions $\cH(\C^n)$. 
  The class $\cE^{\{\ga^s\}}$ coincides with the Gevrey class $\cG^{1+s}$.
\end{examples}

\subsection{Associated sequences}
For $\om \in \sW$ and each $\rh>0$ consider the sequence $\Om^\rh \in \R_{>0}^\N$ defined by 
\[
\Om^\rh_k:=\tfrac{1}{k!}\exp(\tfrac{1}{\rh} \vh^*(\rh k)).
\]
By the properties of $\vh^*$,
each $\Om^\rh$ is weakly log-convex, $(k! \Om^\rh_k)^{\frac{1}{k}} \nearrow \infty$, and  
$\Om^\rh \le \Om^\si$ if $\rh \le \si$. 
By \eqref{eq:cx0},
\begin{equation} \label{eq:cx}
j!\Om^{\rh}_j k!\Om^{\rh}_k \le (j+k)! \Om^\rh_{j+k} \le j!\Om^{2\rh}_j k!\Om^{2\rh}_k, \quad j,k \in \N. 
\end{equation}
In particular, $\Om^\rh_{k+1} \le C \Om^{2\rh}_k$ for all $k$, where $C>0$ is a constant depending on $\rh$.

With $\Om^\rh$ we may associate the function $\om_\rh := \log \o T_{\Om^\rh}$,  
cf.\ \cite[(3.1)]{Komatsu73}; then
\[
\om_\rh(t) = \sup_{k \in \N} (k \log t-\tfrac{1}{\rh}\vh^*(\rh k)) 
\le \sup_{s \ge 0} (s \log t-\tfrac{1}{\rh}\vh^*(\rh s)) = \tfrac{1}{\rh} \om(t).  
\]

\begin{lemma} \label{lem:ass}
  For $\om \in \sW$ we have $\om \approx \om_\rh$ for all $\rh>0$.
\end{lemma}

\begin{demo}{Proof}
  It suffices to show that $\om \approx \om_1$; for arbitrary $\rh>0$ replace $\om$ by $\tfrac{1}{\rh}\om$.
  By \cite[1.8.III]{Mandelbrojt52}, 
  \[
  \om_1(t) = \sup_{k \in \N} (k \log t-\vh^*(k)) = k_t \log t-\vh^*(k_t), 
  \]   
  where $k_t \in \N$ is such that $\vp_{k_t} \le  t < \vp_{k_t+1}$ and $\vp_k := k \Om^1_k/\Om^1_{k-1} \nearrow \infty$. 
  Consider the function $f_t : [0,\infty) \to \R$ given by $f_t(s) = s \log t-\vh^*(s)$, which is concave (for $t\ge 1$) since $\vh^*$ 
  is convex. Concavity of $f_t$ shows that $\om(t) = \sup_{s\ge 0} f_t(s) = f_t(s_t)$ for a point $s_t \in (k_t-1,k_t+1)$. 
  
  Assume that $s_t \in (k_t,k_t+1)$. By concavity of $f_t$ and since $f_t(0)=0$, we find 
  \[
  f_t(s_t) \le \tfrac{f_t(k_t)}{k_t} s_t \le \tfrac{f_t(k_t)}{k_t} (k_t+1) \le 2 f_t(k_t) 
  \]
  and hence $\om(t) \le 2 \om_1(t)$ for sufficiently large $t$. The case $s_t \in (k_t-1,k_t)$ is similar.
\qed\end{demo}

\begin{corollary} \label{cor:omqa}
  For $\om \in \sW$ we have:
  \begin{enumerate}[$(1)$]
    \item $\om$ is quasianalytic if and only if each (equivalently, some) $\Om^\rh$ is quasianalytic.
    \item $\om$ satisfies \thetag{\hyperref[om_6]{$\om_6$}} if and only if each (equivalently, some) $\Om^\rh$ has moderate growth.
  \end{enumerate}
\end{corollary}

\begin{demo}{Proof}
  This follows from Lemma~\ref{lem:ass}, \cite[Lemma~4.1]{Komatsu73}, and \cite[Prop.~3.6]{Komatsu73}.
\qed\end{demo}

\begin{lemma}
  For $\om \in \sW$ we have 
  \begin{gather}
    \forall \si >0 ~\exists H \ge 1 ~\forall \rh>0 ~\exists C\ge 1 ~\forall k \in \N : 
    \si^k \Om^\rh_k \le C \Om^{H\rh}_k.
    \label{eq:abs2}
  \end{gather}
  Moreover, $\om \in \sW$ satisfies \thetag{\hyperref[om_6]{$\om_6$}} if and only if
  \begin{gather}
  \forall \rh>0 ~\forall \ta>0 : \Om^\rh \approx \Om^\ta.  
  \label{eq:abs4}
  \end{gather}
  If $\om \in \sW$ satisfies \thetag{\hyperref[om_8]{$\om_8$}} then
  \begin{gather}
  \exists C>1 ~\forall \rh>0 :  \Om^{\tfrac{\rh}{C}} \lhd \Om^\rh.  
  \label{eq:abs5}
  \end{gather}
\end{lemma}

It follows that \thetag{\hyperref[om_8]{$\om_8$}} is an obstruction for \thetag{\hyperref[om_6]{$\om_6$}}.

\begin{demo}{Proof}
  The following inequality is well-known (see e.g.\ \cite[p.~404]{FernandezGalbis06}):
  \begin{equation}
    \exists L\ge 1 ~\forall t \ge 0 ~\forall s \in \N :
     L^s \vh^*(t)+ s L^s t \le  \vh^*(L^s t) + \sum_{i=1}^s L^i. 
    \label{eq:abs0}
  \end{equation}
  For the reader's convenience we give a short proof.
  By \thetag{\hyperref[om_1]{$\om_1$}}, there exists $L_1\ge1$ such that $\om(2t) \le L_1 \om(t)+L_1$ for all $t\ge0$ and hence
  there exists $L\ge1$ such that $\vh(t+1) \le L \vh(t)+L$ for all $t\ge0$. Thus, for $t \ge 0$,
  \begin{align*}
    \vh^*(Lt) + L = \sup_{s\ge 0} (Lts - (\vh(s)-L)) \ge \sup_{s\ge 1} (Lts - L \vh(s-1)) = L \vh^*(t)+Lt,
  \end{align*}
  and \eqref{eq:abs0} follows by iteration.
  
  By choosing $s$ such that $e^s\ge \si$ and by setting $t:=\rh k$,  
  $H:=L^s$ and $C := \exp(\tfrac{1}{H \rh} \sum_{i=1}^s L^i)$, we see that  
  \eqref{eq:abs0} implies \eqref{eq:abs2}. 

  Let us prove that \thetag{\hyperref[om_6]{$\om_6$}} implies \eqref{eq:abs4}.
  By \thetag{\hyperref[om_6]{$\om_6$}} there exists a constant $H \ge 1$ such that $2 \om(t) \le \om(H t)+ H$ for all $t \ge 0$, and, 
  consequently, as $\om|_{[0,1]}=0$, 
    \begin{align*}
     \vh^*(t) &= \sup_{s\ge 0} (ts -\om(e^s)) = \sup_{s \in \R} (ts -\om(e^s)) = \sup_{u\ge 0} (t \log u - \om(u)) \\
      &\ge \sup_{u\ge 0} (t \log u - \tfrac{1}{2} \om(H u)) - \tfrac{1}{2} H  
      = \tfrac{1}{2} \vh^*(2 t) - t \log H - \tfrac{1}{2} H. 
    \end{align*}
  By setting $t := \rh k$, we may conclude that 
  \[
  \exists H \ge 1 ~\forall \rh>0 ~\forall k \in \N : \Om^{2\rh}_k \le e^{\frac{H}{2\rh}} H^k \Om^\rh_k
  \]  
  which implies $\Om^{2\rh} \preceq \Om^\rh$ for all $\rh>0$. Iteration and the fact that $\Om^\rh \le \Om^{2\rh}$ yield 
  $\Om^{2^n\rh} \approx \Om^\rh$ for all $\rh>0$ and all $n \in \N$, and \eqref{eq:abs4} follows.
  
  Conversely assume \eqref{eq:abs4} which means that
  \[
  \forall \rh>0 ~\forall \ta>0 ~\exists C>0 ~\forall k \in \N : \tfrac{1}{\rh} \vh^*(\rh k) \le C k + \tfrac{1}{\ta} \vh^*(\ta k).  
  \]   
  By \eqref{eq:cx0}, we may conclude that
  \begin{align*}
  \forall \rh>0 ~\forall \ta>0 ~\exists D>0 ~\forall t \ge 0 : \tfrac{1}{\rh} \vh^*(\rh t) \le D t + D + \tfrac{1}{2\ta} \vh^*(2\ta t).
  \end{align*} 
  Thus 
  \begin{align*}
    \tfrac{1}{2 \ta} \vh(t) &= \sup_{s \ge 0} (ts - \tfrac{1}{2 \ta} \vh^*(2 \ta s)) \le 
    \sup_{s \ge 0} (ts +Ds - \tfrac{1}{\rh} \vh^*(\rh s)) + D 
    = \tfrac{1}{\rh} \vh(t+D) +D,
  \end{align*}
  and, hence,
  \[
  \tfrac{1}{2 \ta} \om(t) \le \tfrac{1}{\rh} \om(e^D t) +D.
  \] 
  Setting $\rh = 4$ and $\ta=1$ implies \thetag{\hyperref[om_6]{$\om_6$}}.
  
  Let us prove \eqref{eq:abs5}. By \thetag{\hyperref[om_8]{$\om_8$}} there exist constants $C,H>0$ such that
  \begin{align*}
    C \vh^*(\tfrac{2t}{C}) &= \sup_{u\ge 0} (2 t \log u - C\om(u)) = \sup_{u\ge 0} (2 t \log u - C\om(H u)) + 2 t \log H \\
    &\le \sup_{u\ge 0} (2 t \log u - \om(u^2)) + 2 t \log H +C 
    = \vh^*(t)  + 2 t \log H +C.
  \end{align*}
  By setting $t := \rh k$ we find that for all $\rh>0$ and all $k \in \N$
  \[
  (2k)! \Om^{\frac{\rh}{C}}_{2k} \le  e^{\tfrac{C}{\rh}} H^{2k} k! \Om^{\rh}_{k}.
  \]
  Thus the sequence $L = (L_k)$ defined by $k! L_k := (2k)! \Om^{\tfrac{\rh}{C}}_{2k} \ge (k! \Om^{\tfrac{\rh}{C}}_k)^2$ satisfies
  $\Om^{\tfrac{\rh}{C}} \lhd L \preceq \Om^\rh$, which implies \eqref{eq:abs5}.
\qed\end{demo}

\begin{theorem} \label{thm:rep}
  Let $\om \in \sW$, 
  let $U \subseteq \R^n$ be open, and let $K \subseteq U$ be compact. Then:
  \begin{enumerate}[$(1)$]
    \item For each $\rh >0$ we have $\cE^{\{\Om^\rh\}}(U) \subseteq \cE^{\{\om\}}(U)$ and $\cE^{(\om)}(U) \subseteq \cE^{(\Om^\rh)}(U)$ 
    with continuous inclusion.
    \item We have as locally convex spaces
  \[
  \cE^{(\om)}(U) = \varprojlim_{\rh>0} \cE^{(\Om^\rh)}(U) \quad 
  \text{and} \quad \cE^{\{\om\}}(K) = \varinjlim_{\rh>0} \cE^{\{\Om^\rh\}}(K).
  \]  
    \item[\thetag{3}] $\om$ satisfies \thetag{\hyperref[om_6]{$\om_6$}} if and only if 
    $\cE^{[\Om^\rh]}(U) = \cE^{[\om]}(U)$, for each $\rh >0$, as locally convex spaces.
    \item[\thetag{4}] If $\om$ satisfies \thetag{\hyperref[om_8]{$\om_8$}}, then also
    \begin{align*}
    \cE^{(\om)}(U) &= \varprojlim_{\rh>0} \cE^{(\Om^\rh)}(U) = \varprojlim_{\rh>0} \cE^{\{\Om^\rh\}}(U) \quad 
    \text{and} \quad \\ 
    \cE^{\{\om\}}(K) &= \varinjlim_{\rh>0} \cE^{\{\Om^\rh\}}(K) = \varinjlim_{\rh>0} \cE^{(\Om^\rh)}(K)
    \end{align*}
    as locally convex spaces.
  \end{enumerate}
\end{theorem}

\begin{demo}{Proof}
  \thetag{1} Let $\rh>0$ be fixed.
  If $f \in \cE^{\{\Om^\rh\}}(U)$ then for each compact $K \subseteq U$ there exists $\si>0$ such that $\|f\|^{\Om^\rh}_{K,\si} < \infty$.
  By \eqref{eq:abs2}, there exist constants $H, C \ge 1$ such that
  \[
  \infty >
  C \|f\|^{\Om^\rh}_{K,\si} \ge 
  \|f\|^{\Om^{H\rh}}_{K,1} = \|f\|^{\om}_{K,H\rh},
  \]
  whence $f \in \cE^{\{\om\}}(U)$.
  
  Assume that $f \in \cE^{(\om)}(U)$. Let $\rh>0$ and $\si >0$ be fixed. By \eqref{eq:abs2}, there exist constants $H, C \ge 1$ such that 
  $\Om^\rh_k \le C \si^k \Om^{H\rh}_k$ for all $k$. Since $f \in \cE^{(\om)}(U)$, for each compact $K \subseteq U$ we have 
  $\|f\|^{\om}_{K,\frac{\rh}{H}} < \infty$, and, thus,
  \[
  \infty > C \|f\|^{\om}_{K,\frac{\rh}{H}} = C \|f\|^{\Om^{\frac{\rh}{H}}}_{K,1} \ge \|f\|^{\Om^\rh}_{K,\si}. 
  \]   
  Since $\si>0$ was arbitrary, we may conclude that $f \in \cE^{(\Om^\rh)}(U)$.
  
  \thetag{2} follows from \thetag{1},
  since the inclusions $\cE^{(\om)}(U) \supseteq \varprojlim_{\rh>0} \cE^{(\Om^\rh)}(U)$ and
  $\cE^{\{\om\}}(K) \subseteq \varinjlim_{\rh>0} \cE^{\{\Om^\rh\}}(K)$ are clear and continuous by definition.
  
  \thetag{3} follows from \thetag{2}, \eqref{eq:abs4}, and Proposition~\ref{prop:union}\thetag{1}.

  \thetag{4} is a direct consequence of \thetag{2}, \eqref{eq:abs5}, and Proposition~\ref{prop:union}\thetag{1}. 
\qed\end{demo}

\begin{corollary} \label{cor:rep}
  Let $\om \in \sW$ and let $U \subseteq \R^n$ be open. 
  Then $\cE^{[\om]}(U) = \cE^{[\fW]}(U)$ as locally convex spaces, 
  where the weight matrix $\fW := \{\Om^\rh : \rh>0\}$ 
  satisfies
  \begin{itemize}
    \item \thetag{\hyperref[fM_(mg)]{$\fM_{(\on{mg})}$}} and \thetag{\hyperref[fM_{mg}]{$\fM_{\{\on{mg}\}}$}},
    \item \thetag{\hyperref[fM_(alg)]{$\fM_{(\on{alg})}$}} and \thetag{\hyperref[fM_{alg}]{$\fM_{\{\on{alg}\}}$}},
    \item \thetag{\hyperref[fM_(L)]{$\fM_{(\on{L})}$}} and \thetag{\hyperref[fM_{L}]{$\fM_{\{\on{L}\}}$}}.
  \end{itemize}
  If $\om$ satisfies \thetag{\hyperref[om_4]{$\om_4$}}, respectively \thetag{\hyperref[om_5]{$\om_5$}}, 
  then $\fW$ satisfies \thetag{\hyperref[fM_H]{$\fM_{\cH}$}}, respectively 
  \thetag{\hyperref[fM_(Com)]{$\fM_{(C^\om)}$}}.
  If $\om$ satisfies \thetag{\hyperref[om_4]{$\om_8$}}, then $\fW$ satisfies \thetag{\hyperref[fM_(BR)]{$\fM_{(\on{BR})}$}} 
  and \thetag{\hyperref[fM_{BR}]{$\fM_{\{\on{BR}\}}$}}.    
\end{corollary}

\begin{demo}{Proof}
  This is an immediate consequence of Theorem~\ref{thm:rep}, \eqref{eq:cx}, and \eqref{eq:abs2}.

  For $\om(t)=\max\{0,t-1\} \approx t$ we find $\vh^*(t)=t \log t-t+1$, for $t\ge 1$, $\vh^*|_{[0,1]}=0$, 
  and it is easy to see that \thetag{\hyperref[om_4]{$\om_4$}} implies 
  \thetag{\hyperref[fM_H]{$\fM_{\cH}$}} and 
  \thetag{\hyperref[om_5]{$\om_5$}} implies \thetag{\hyperref[fM_(Com)]{$\fM_{(C^\om)}$}}, 
  by Lemma~\ref{lem:equi}.  
  That \thetag{\hyperref[om_4]{$\om_8$}} implies \thetag{\hyperref[fM_(BR)]{$\fM_{(\on{BR})}$}} 
  and \thetag{\hyperref[fM_{BR}]{$\fM_{\{\on{BR}\}}$}} follows from \eqref{eq:abs5}.
\qed\end{demo}

\begin{lemma} \label{lem:equi}
For $\om, \si \in \sW$ we have: 
\begin{enumerate}[$(1)$]
  \item If $\om \preceq \si$ 
  then $\exists H \ge 1 ~\forall \rh>0 ~\exists C>0 : 
  \Om^\rh \le C \Si^{H \rh}$.
  \item If $\om \lhd \si$ 
  then $\forall H >0 ~\forall \rh>0 ~\exists C>0 : 
  \Om^\rh \le C \Si^{H \rh}$.
\end{enumerate}
Here $\Si^\rh$ are the sequences associated with $\si$.
\end{lemma}

\begin{demo}{Proof}
  \thetag{1} If $\om \preceq \si$ then there exists $H \ge 1$ such that $\si(t)\le H \om(t)+H$ for all $t\ge 0$, and thus also 
  $\vh_\si(t)\le H \vh_\om(t)+H$ and finally $H \vh_\om^*(t) \le \vh_\si^*(Ht) + H$. Setting $t= \rh k$ gives the assertion.
  
  \thetag{2} If $\om \lhd \si$ then for all $H>0$ there exists $D>0$ such that $\si(t) \le H \om(t) + D$ for all $t\ge 0$, and thus 
  $H \vh_\om^*(t) \le \vh_\si^*(Ht) + D$ as in \thetag{1}. 
  Setting $t= \rh k$ gives the assertion.
\qed\end{demo}

\begin{corollary} 
  For $\om, \si \in \sW$ we have:
  \begin{enumerate}[$(1)$]
    \item  $\om \preceq \si \Rightarrow \cE^{[\om]} \subseteq \cE^{[\si]}$ and 
    $\cE^{[\om]}(\R) \subseteq \cE^{[\si]}(\R) \Rightarrow \om \preceq \si$. 
    \item  $\om \lhd \si \Rightarrow \cE^{\{\om\}} \subseteq \cE^{(\si)}$ and 
    $\cE^{\{\om\}}(\R) \subseteq \cE^{(\si)}(\R) \Rightarrow \om \lhd \si$.
    \item  There is no $\si \in \sW$ such that 
    $\cE^{(\om)}(\R) \subsetneq \cE^{[\si]}(\R) \subsetneq \cE^{\{\om\}}(\R)$.
  \end{enumerate}
\end{corollary}

\begin{demo}{Proof}
  \thetag{1}
  If $\fW := \{\Om^\rh : \rh>0\}$ and $\fS := \{\Si^\rh : \rh>0\}$, where $\Si^\rh$ are the sequences associated with $\si$, then 
  in view of Proposition~\ref{prop:fMincl} and Corollary~\ref{cor:rep} 
  it suffices to show
  \begin{enumerate}[$(1')$]
    \item  $\om \preceq \si$ if and only if $\fW [\preceq] \fS$.  
  \end{enumerate}
  If $\om \preceq \si$ then Lemma~\ref{lem:equi} implies $\fW (\preceq) \fS$ as well as $\fW \{\preceq\} \fS$.
  
  Conversely, assume $\fW \{\preceq\} \fS$, i.e., using \eqref{eq:abs2}, 
  \[
  \forall \rh>0 ~\exists \ta>0 ~\exists C>0 ~\forall k \in \N : \tfrac{1}{\rh} \vh^*_\om(\rh k) 
  \le \tfrac{1}{\ta} \vh^*_\si(\ta k) + C, 
  \]
  and, by \eqref{eq:cx0},
  \[
  \forall \rh>0 ~\exists \ta>0 ~\exists D>0 ~\forall t\ge 0 : \tfrac{1}{\rh} \vh^*_\om(\rh t) 
  \le \tfrac{1}{2\ta} \vh^*_\si(2\ta t) + D. 
  \]
  Thus 
  \begin{align*}
    \tfrac{1}{2 \ta} \vh_\si(t) &= \sup_{s \ge 0} (ts - \tfrac{1}{2 \ta} \vh^*_\si(2 \ta s)) \le 
    \sup_{s \ge 0} (ts - \tfrac{1}{\rh} \vh^*_\om(\rh s)) + D 
    = \tfrac{1}{\rh} \vh_\om(t) +D,
  \end{align*}
  and, hence,
  \begin{equation} \label{eq:Oo}
      \tfrac{1}{2 \ta} \si(t) \le \tfrac{1}{\rh} \om(t) +D,
  \end{equation} 
  which implies $\si(t) = O(\om(t))$ as $t \to \infty$, i.e., $\om \preceq \si$.
  
  If $\fW (\preceq) \fS$, then the same arguments yield \eqref{eq:Oo}, but with swapped quantifiers:
  \[
  \forall \ta>0 ~\exists \rh>0 ~\exists D>0 ~\forall t\ge0 : \tfrac{1}{2 \ta} \si(t) \le \tfrac{1}{\rh} \om(t) +D.
  \]
  Again this implies $\om \preceq \si$.
  
  \thetag{2}
  If $\om \lhd \si$ then Lemma~\ref{lem:equi} implies $\cE^{\{\om\}} \subseteq \cE^{(\si)}$. 
  Conversely, if $\cE^{\{\om\}}(\R) \subseteq \cE^{(\si)}(\R)$, then
  $\cE^{\{\Om^\rh\}}$ admits a characteristic function and is contained in $\cE^{(\si)}$, thus 
  \[
  \forall \rh>0 ~\forall \ta>0 ~\exists C>0 ~\forall k \in \N : \tfrac{1}{\rh} \vh^*_\om(\rh k) 
  \le \tfrac{1}{\ta} \vh^*_\si(\ta k) + C.
  \]
  As in \thetag{1} we may derive that for all $\rh,\ta > 0$ there is $D>0$ such that \eqref{eq:Oo} for all $t\ge 0$, 
  hence $\si(t) = o(\om(t))$ as $t \to \infty$, i.e., $\om \lhd \si$.

  \thetag{3} If $\cE^{(\om)}(\R) \subseteq \cE^{\{\si\}}(\R)$, then $\fW (\preceq\} \fS$, 
  by Corollary~\ref{cor:rep} and Proposition~\ref{prop:fMincl}.
  Similarly as in \thetag{1} we may derive that there exist $\rh,\ta>0$ such that \eqref{eq:Oo}, 
  and so $\om \preceq \si$. This and \thetag{1} imply the assertion.
\qed\end{demo}

As $\cE^{\{t\}}(U)=C^\om(U)$ and $\cE^{(t)}(U) = \cH(\C^n)$ (via restriction), condition \thetag{\hyperref[om_4]{$\om_4$}} is equivalent to  
$C^\om \subseteq \cE^{\{\om\}}$ and condition \thetag{\hyperref[om_5]{$\om_5$}} is equivalent to $C^\om \subseteq \cE^{(\om)}$.

\subsection{Intersection and union of all non-quasianalytic Gevrey classes} \label{ssec:Gevrey}
For the weight matrix $\fG = \{G^s : s > 0\}$ with $G^s=(G^s_k) = ((k!)^s)$
\begin{equation} \label{eq:Gint}
    \cE^{(\fG)}(U) = \bigcap_{s>0} \cG^{1+s}(U), 
    \quad U \subseteq \R^n \text{ open,}
  \end{equation}
  is the intersection and 
  \begin{equation} \label{eq:Gun}
    \cE^{\{\fG\}}(K) = \bigcup_{s>0} \cG^{1+s}(K), \quad K \subseteq \R^n \text{ compact,}  
\end{equation}
is the union of all 
non-quasianalytic Gevrey classes $\cG^{1+s} = \cE^{\{G^s\}}$
(as locally convex spaces). 
Indeed $G^s \lhd G^{s'}$ for all $s < s'$ (so $\fG$ satisfies \thetag{\hyperref[fM_(BR)]{$\fM_{(\on{BR})}$}} 
and \thetag{\hyperref[fM_{BR}]{$\fM_{\{\on{BR}\}}$}}), and hence we get \eqref{eq:Gint}
\[
\cE^{(\fG)}(U) = \bigcap_{s>0} \cE^{(G^s)}(U) = \bigcap_{s>0} \cE^{\{G^s\}}(U) = \bigcap_{s>0} \cG^{1+s}(U), 
\]
while \eqref{eq:Gun} is evident by definition. 
Note that $\cE^{(\fG)}$, and hence also $\cE^{\{\fG\}}$, 
is non-quasianalytic; in fact, the sequence $L=(L_k)$ defined by $k! L_k := k^k(\log(k+e))^{2k}$ is non-quasianalytic 
and satisfies $L \lhd G^s$ for all $s>0$, and, as $(k! L_k)^{\frac{1}{k}}$ is increasing, $\cE^{[L]}$ is non-quasianalytic, 
by the Denjoy--Carleman theorem.  

The following theorem shows that 
there exist spaces $\cE^{[\fM]}$ that are different from $\cE^{[M]}$ as well as from $\cE^{[\om]}$.

\begin{theorem} \label{thm:ext}
  Neither $\cE^{(\fG)}(\R)$ nor $\cE^{\{\fG\}}(\R)$ coincides (as vector space) with 
  $\cE^{(M)}(\R)$, $\cE^{\{M\}}(\R)$, $\cE^{(\om)}(\R)$, or $\cE^{\{\om\}}(\R)$ for any weight sequence $M$ 
  or weight function $\om$. 
\end{theorem}

\begin{demo}{Proof}
  We show first that, given a weight matrix $\fM =\{M^\la : \la \in \La\}$ with $M^\la \not\approx M^\mu$ for all $\la \ne \mu$, 
  there cannot exist a weakly log-convex $M \in \R_{>0}^{\N}$ such that $\cE^{[\fM]}(\R)=\cE^{[M]}(\R)$. 
  Indeed, if there is such $M$, Proposition~\ref{prop:fMincl} implies $M \approx M^\la$ for some $\la$.
  Then, by Proposition~\ref{prop:union}\thetag{1}, 
  \[
    \cE^{(M)}(\R) = \cE^{(\fM)}(\R) = \bigcap_{\la} \cE^{(M^\la)}(\R) \subsetneq  \cE^{(M)}(\R), 
  \]
  and, for compact $K \subseteq \R$, 
  \begin{align*}
    \cE^{\{M\}}(\R) = \cE^{\{\fM\}}(\R) &= \bigcap_K \cE^{\{\fM\}}(K) = \bigcap_K \bigcup_{\la} \cE^{\{M^\la\}}(K) \\
    &\supseteq \bigcup_{\la} \bigcap_K \cE^{\{M^\la\}}(K) = \bigcup_{\la} \cE^{\{M^\la\}}(\R) \supsetneq \cE^{\{M\}}(\R),  
  \end{align*}
  which contradicts the assumption in both cases.

  As $\cE^{(\fG)}(\R)$ contains $C^\om(\R)$ it cannot coincide with $\cE^{(M)}(\R)$ for any weight sequence $M$, 
  by Theorem~\ref{thm:reg} and the first paragraph; neither can $\cE^{\{\fG\}}(\R)$ coincide with $\cE^{\{M\}}(\R)$.     

  If there exists $\om \in \sW$ such that $\cE^{(\fG)}(\R)=\cE^{(\om)}(\R)$, 
  then Proposition~\ref{prop:fMincl} implies that  
  for each $\rh>0$ there exist $s,\rh'>0$ such that 
  \begin{equation} \label{eq:OGO}
    \Om^{\rh'} \preceq G^{s} \preceq \Om^{\rh},
  \end{equation}
  and thus, by Proposition~\ref{prop:union}\thetag{1},
  \[
    \cE^{\{\Om^{\rh'}\}} \subseteq \cG^{1+s} \subseteq \cE^{\{\Om^{\rh}\}}.
  \]
  Since $\cG^{1+s} = \cE^{\{\ga\}}$ with $\ga(t)= t^{\tfrac{1}{1+s}}$,
  using the fact that there exist characteristic $\cE^{\{\Om^\rh\}}$- and $\cE^{\{\Ga^{\ta}\}}$-functions 
  (where $\Ga^\ta$ are the 
  sequences associated with $\ga$), and by \eqref{eq:abs2}, we may conclude that, for all $k$,
  \[
    \tfrac{1}{\rh'} \vh_\om^*(\rh' k) \le \tfrac{1}{\ta} \vh_\ga^*(\ta k) + C \quad \text{ and } \quad
    \tfrac{1}{\ta} \vh_\ga^*(\ta k) \le \tfrac{1}{H\rh} \vh_\om^*(H\rh k) + D,
  \]
  for suitable constants $\ta,C,D,H$.
  As in the derivation of \eqref{eq:Oo}
  this implies $\om \approx \ga$ and hence $\cE^{(\fG)}(\R)=\cE^{(\om)}(\R) = \cE^{(\ga)}(\R) = \cE^{(G^s)}(\R)$, 
  a contradiction. Thus there is no $\om \in \sW$ with $\cE^{(\fG)}(\R)=\cE^{(\om)}(\R)$.

  If there exists $\om \in \sW$ such that $\cE^{\{\fG\}}(\R)=\cE^{\{\om\}}(\R)$, then Proposition~\ref{prop:fMincl} implies that  
  for each $\rh'>0$ there exist $s,\rh>0$ such that \eqref{eq:OGO}. Then the same arguments show 
  $\om \approx \ga$ and hence $\cE^{\{\fG\}}(\R)=\cE^{\{\om\}}(\R) = \cE^{\{\ga\}}(\R) = \cG^{1+s}(\R)$, 
  a contradiction. Thus there is no $\om \in \sW$ with $\cE^{\{\fG\}}(\R)=\cE^{\{\om\}}(\R)$.

  For the remaining cases note that   
  $\fM (\preceq\} \fN \{\lhd) \fM$ as well as $\fM \{\lhd) \fN (\preceq\} \fM$ is impossible
  for any two weight matrices $\fM,\fN \in \sM$. 
  This fact together with Proposition~\ref{prop:fMincl} (and Theorem~\ref{thm:reg}) 
  implies that there is no weight sequence $M$ and 
  no weight function $\om$ so that 
  $\cE^{(\fG)}(\R)=\cE^{\{M\}}(\R)$, $\cE^{(\fG)}(\R)=\cE^{\{\om\}}(\R)$, 
  $\cE^{\{\fG\}}(\R)=\cE^{(M)}(\R)$, or $\cE^{\{\fG\}}(\R)=\cE^{(\om)}(\R)$.
  The proof is complete.
\qed\end{demo}

\begin{corollary}
  Composition is continuous on the intersection of all non-quasianalytic Gevrey classes.
  More precisely,
  $\on{comp}^{(\fG)}$ is continuous, $\cE^{\{\fG\}}(\R^p,f)$, for $f \in \cE^{\{\fG\}}(\R^q,\R^r)$, is continuous, and $\on{comp}^{\{\fG\}}$ is sequentially continuous.  
\end{corollary}

\begin{demo}{Proof}
  This follows from Theorem~\ref{thm:compC0} and Theorem~\ref{thm:ext}.
\qed\end{demo}

We expect that $\on{comp}^{[\fG]}$ is even $\cE^{[\fG]}$, see
Remark~\ref{rem:fc}.

\begin{remark}
  More autonomous spaces $\cE^{[\fM]}$ can be produced by choosing the weight matrix $\fM := \{M^\la : \la>0\}$ such that 
  each $M^\la$ has moderate growth, satisfies $\varliminf (M^\la_k)^{\frac{1}{k}}>0$ and 
  $\varliminf \mu^\la_{nk}/\mu^\la_k>1$ for some $n \in \N$ with $\mu^\la_k=kM^\la_k/M^\la_{k-1}$, 
  and $M^\la \not \approx M^\mu$ for $\la \ne \mu$. 
  Here we may use the comparison theorems in \cite{BMM07} and argue as above. 
\end{remark}

\section{Stability under composition of \texorpdfstring{$\cE^{[\om]}$}{E[om]}} \label{sec:wfc}

Stability under composition of $\cE^{[\om]}$ was characterized in \cite{FernandezGalbis06} for non-quasianalytic weights $\om$.
In this section we apply the characterization obtained by means of the associated weight matrix $\fW = \{\Om^\rh : \rh>0\}$ and relate it 
to the results of \cite{FernandezGalbis06}.

\begin{lemma} \label{lem:sa}
  If $\om \in \sW$ is sub-additive, 
  then for each $\rh>0$ we have $(\Om^\rh)^\o \preceq \Om^{2\rh}$.
\end{lemma}

Then the weight matrix $\fW$ satisfies \thetag{\hyperref[fM_(FdB)]{$\fM_{(\on{FdB})}$}} and 
\thetag{\hyperref[fM_{FdB}]{$\fM_{\{\on{FdB}\}}$}}.

\begin{demo}{Proof}
  Sub-additivity of $\om$ implies
  \begin{equation} \label{eq:sa}
  \Om^\rh_j \Om^\rh_k \le \Om^\rh_{j+k}, \quad j,k \in \N,
  \end{equation}
  see\ \cite[Lemma~3.3]{FernandezGalbisJornet04}.
  Indeed, we have 
  $\exp(\tfrac{1}{\rh}\vh^*(\rh k)) = \sup_{s\ge 1} s^k \exp(-\tfrac{1}{\rh}\om(s))$
  and hence, using sub-additivity of $\om$,
  \begin{align*}
  \Om^\rh_j \Om^\rh_k 
  \le \sup_{s,t\ge 1}\frac{s^jt^k}{j! k!} \exp(-\tfrac{1}{\rh}\om(s+t))
  \le \sup_{s,t\ge 1}\frac{(s+t)^{j+k}}{(j+k)!} \exp(-\tfrac{1}{\rh}\om(s+t))
  \le \Om^\rh_{j+k}.
  \end{align*}
    
  By \eqref{eq:cx}, \eqref{eq:sa} and since $\Om^\rh \le \Om^{2\rh}$, we get, for $\al_i \in \N_{>0}$ with $\al_1+\cdots+ \al_j=k$,
  \begin{align*}
    \Om^\rh_j \Om^\rh_{\al_1} \cdots \Om^\rh_{\al_j} \le C^j \Om^{2\rh}_j \Om^{2\rh}_{\al_1-1} \cdots \Om^{2\rh}_{\al_j-1} 
    \le C^j \Om^{2\rh}_k  
  \end{align*}
  which implies the assertion.
\qed\end{demo}

\begin{theorem} \label{thm:compomR}
For $\om \in \sW$ satisfying \thetag{\hyperref[om_4]{$\om_4$}} the following are equivalent:
\begin{enumerate}[$(1)$]
  \item $\cE^{\{\om\}}$ is stable under composition.
  \item For each $\rh>0$ there is $\ta>0$ so that $(\Om^\rh)^\o \preceq \Om^\ta$, i.e.,  
  $\fW$ satisfies \thetag{\hyperref[fM_{FdB}]{$\fM_{\{\on{FdB}\}}$}}.
  \item There exists a sub-additive $\tilde \om \in \sW$ so that $\om \approx \tilde \om$.
  \item $\om$ satisfies \thetag{\hyperref[om_7]{$\om_7$}}.
\end{enumerate}
\end{theorem}

\begin{demo}{Proof}
  $\thetag{1} \Leftrightarrow \thetag{2}$ follows from Theorem~\ref{thm:{comp}} and Corollary~\ref{cor:rep}.
  
  $\thetag{3} \Leftrightarrow \thetag{4}$ See \cite[Prop.~1.1]{PetzscheVogt84} and \cite[Lemma~1]{Peetre70}.
  
  $\thetag{3} \Rightarrow \thetag{2}$ follows from Lemma~\ref{lem:sa}. 
  
  $\thetag{2} \Rightarrow \thetag{3}$ 
  The proof is inspired by \cite[Prop.~2.3]{FernandezGalbis06} which treats the non-quasianalytic case.
  We do not assume non-quasianalyticity (or quasianalyticity) and use Claim \ref{cl:1}
  to remedy the lack of $\cE^{\{\om\}}$-functions of compact support.
  If $\om$ does not satisfy \thetag{\hyperref[om_7]{$\om_7$}}, then there exist increasing sequences $(k_n) \in \N^\N$ and 
  $(t_n) \in \R_{>0}^\N$ so that 
  \begin{equation} \label{eq:FG}
    \om(k_n t_n)  \ge n^2 k_n \om(t_n).
  \end{equation}  
  Set $a_n := e^{-n \om(t_n)}$ and $f_n(x) := a_n e^{i t_n x}$, $x \in \R$. Then
  \begin{align*}
    \|f_n\|^\om_{\R,\rh} &= a_n \sup_{j \in \N} (t_n^j \exp(-\tfrac{1}{\rh} \vh^*(\rh j))) = 
    a_n  \exp \sup_{j \in \N} (j \log t_n -\tfrac{1}{\rh} \vh^*(\rh j)) \\
    &= e^{-n \om(t_n)} e^{\om_\rh(t_n)} \le e^{-(n- \frac{1}{\rh})\om(t_n)}
  \end{align*}
  and so $\{f_n : n \in \N\}$ is bounded in $\cE^{\{\om\}}(\R,\C)$ (even in $\cE^{(\om)}(\R,\C)$). 
  The set $\{\C \ni z \mapsto z^k : k \in \N\}$ forms a bounded subset of $\cE^{\{\om\}}(\mathbb D,\C)$, where $\mathbb D \subseteq \C$ is the unit 
  disk and where
  we identify $\C \cong \R^2$). 
  Indeed, for $|z|\le r<1$ choose $\rh>0$ such that $r+\frac{1}{\rh}<1$, and thus
  \[
  \sup_{j \in \N} \frac{|\p_z^j z^k|}{\rh^j j!} \le \sup_{j \le  k} \binom{k}{j} r^{k-j} \frac{1}{\rh^j} \le \Big(r+\frac{1}{\rh}\Big)^k.
  \]
  So $\{z \mapsto z^k : k \in \N\}$ is bounded in $C^\om(\mathbb D,\C)$ and, by \thetag{\hyperref[om_4]{$\om_4$}}, 
  in $\cE^{\{\om\}}(\mathbb D,\C)$.
  Since $\fW$ satisfies \thetag{\hyperref[fM_{FdB}]{$\fM_{\{\on{FdB}\}}$}} by assumption \thetag{2},
  we may conclude, from Claim \ref{cl:1}, that the set $\{f_n^k : n,k \in \N\}$ is bounded in $\cE^{\{\om\}}(\R,\C)$. 
  Thus there exists $\rh>0$ such that 
  \begin{align*}
   \infty &> \sup_{n,k,j \in \N} |(f_n^k)^{(j)}(0)| \exp(-\tfrac{1}{\rh} \vh^*(\rh j)) 
   = \sup_{n,k,j \in \N} a_n^k (t_n k)^j \exp(-\tfrac{1}{\rh} \vh^*(\rh j)) \\
   &= \sup_{n,k \in \N} a_n^k e^{\om_\rh (t_n k)} \ge D \sup_{n,k \in \N} a_n^k e^{\frac{1}{C}\om (t_n k)} 
   = D \sup_{n,k \in \N} e^{-nk \om(t_n)+ \frac{1}{C}\om (t_n k)},
  \end{align*}
  for constants $C,D>0$, 
  by Lemma \ref{lem:ass},
  which contradicts \eqref{eq:FG}.
\qed\end{demo}

\begin{theorem} \label{thm:compomB}
For $\om \in \sW$ satisfying \thetag{\hyperref[om_4]{$\om_4$}} the following are equivalent:
\begin{enumerate}[$(1)$]
  \item $\cE^{(\om)}$ is stable under composition.
  \item $\cE^{(\om)}$ is holomorphically closed.
  \item For each $\rh>0$ there exists $\ta>0$ so that $(\Om^\ta)^\o \preceq \Om^{\rh}$, i.e.,  
  $\fW$ satisfies \thetag{\hyperref[fM_(FdB)]{$\fM_{(\on{FdB})}$}}.
  \item There exists $H\ge 1$ so that for each $\rh>0$ we have $(\Om^\rh)^\o \preceq \Om^{H \rh}$. 
  \item There exists a sub-additive $\tilde \om \in \sW$ so that $\om \approx \tilde \om$.
  \item $\om$ satisfies \thetag{\hyperref[om_7]{$\om_7$}}.
\end{enumerate}
\end{theorem}

Note that \thetag{\hyperref[om_4]{$\om_4$}} is needed only for $\thetag{1} \Rightarrow \thetag{2}$.

\begin{demo}{Proof}
  $\thetag{1} \Leftrightarrow \thetag{2} \Leftrightarrow \thetag{3}$ follows from
  Theorem~\ref{thm:fMB} and Corollary~\ref{cor:rep}.
  
  $\thetag{2} \Rightarrow \thetag{6}$ follows from an argument due to \cite{Bruna80/81}, see \cite[p.~405]{FernandezGalbis06}.
  
  $\thetag{5} \Leftrightarrow \thetag{6}$ See \cite[Prop.~1.1]{PetzscheVogt84} and \cite[Lemma~1]{Peetre70}.
    
  $\thetag{5} \Rightarrow \thetag{4}$ follows from Lemma~\ref{lem:sa} and Lemma~\ref{lem:equi}.
  
  $\thetag{4} \Rightarrow \thetag{3}$ is evident.
\qed\end{demo}

\begin{corollary}
  For $\om \in \sW$ satisfying \thetag{\hyperref[om_4]{$\om_4$}} the following are equivalent:
  \begin{enumerate}[$(1)$]
    \item For each $\rh>0$ there exists $\ta>0$ such that $(\Om^\rh)^\o \preceq \Om^\ta$.
    \item For each $\rh>0$ there exists $\ta>0$ such that $(\Om^\ta)^\o \preceq \Om^{\rh}$.
    \item There exists $H\ge 1$ so that for each $\rh>0$ we have $(\Om^\rh)^\o \preceq \Om^{H \rh}$.
  \end{enumerate}
\end{corollary}

\begin{demo}{Proof}
  Combine Theorem~\ref{thm:compomR} and Theorem~\ref{thm:compomB}.
\qed\end{demo}

Special cases of Theorem~\ref{thm:compC0} were proven in 
\cite[4.2 and 4.4]{FernandezGalbis06}:

\begin{corollary}
  Let $\om \in \sW$ satisfy \thetag{\hyperref[om_7]{$\om_7$}}.
  Then $\on{comp}^{(\om)}$ is continuous, $\cE^{\{\om\}}(\R^p,f)$, for $f \in \cE^{\{\om\}}(\R^q,\R^r)$, is continuous, and $\on{comp}^{\{\om\}}$ is sequentially continuous.
\end{corollary}

\begin{demo}{Proof}
  This is a special case of Theorem~\ref{thm:compC0}, by Corollary~\ref{cor:rep}, Theorem~\ref{thm:compomR}, and Theorem~\ref{thm:compomB}. 
\qed\end{demo}

We expect that the mapping $\on{comp}^{[\om]}$ is even $\cE^{[\om]}$, see
Remark~\ref{rem:fc}.

\def\cprime{$'$}
\providecommand{\bysame}{\leavevmode\hbox to3em{\hrulefill}\thinspace}
\providecommand{\MR}{\relax\ifhmode\unskip\space\fi MR }
\providecommand{\MRhref}[2]{%
  \href{http://www.ams.org/mathscinet-getitem?mr=#1}{#2}
}
\providecommand{\href}[2]{#2}


\end{document}